\documentclass[review, compress, 3p]{elsarticle}
\usepackage{lineno,hyperref}
\usepackage{amsfonts}
\usepackage{amssymb}
\usepackage{amsmath}
\usepackage{rotating,float}
\usepackage{algorithm,algorithmic}
\usepackage{graphicx}
\usepackage{rotating}
\usepackage{float}
\usepackage{booktabs}
\usepackage{multirow}
\usepackage{lscape}
\usepackage{multicol}
\usepackage{subfigure}
\usepackage{color}
\usepackage{setspace}
\usepackage{natbib}
\usepackage{bm}
\usepackage{doi}

\makeatletter
\@addtoreset{equation}{section}
\makeatother

\newtheorem{theorem}{Theorem}[section]
\newtheorem{lemma}{Lemma}[section]
\newtheorem{corollary}{Corollary}

\begin{document}

\begin{frontmatter}

\title{A preconditioning technique for all-at-once system from the nonlinear tempered fractional diffusion equation}

\author[address1]{Yong-Liang Zhao}
\ead{uestc\_ylzhao@sina.com}

\author[address1]{Pei-Yong Zhu}
\ead{zpy6940@uestc.edu.cn}

\author[address2]{Xian-Ming Gu}
\ead{guxianming@live.cn}

\author[address1]{Xi-Le Zhao}
\ead{xlzhao122003@163.com}

\author[address1]{Huan-Yan Jian}
\ead{uestc\_hyjian@sina.com}

\address[address1]{School of Mathematical Sciences, \\
University of Electronic Science and Technology of China,
Chengdu, Sichuan 611731, P.R. China}
\address[address2]{School of Economic Mathematics/Institute of Mathematics, \\
Southwestern University of Finance and Economics, Chengdu, Sichuan 611130, P.R. China}

\begin{abstract}
An all-at-once linear system arising from the nonlinear tempered fractional
diffusion equation with variable coefficients is studied. Firstly, the nonlinear
and linearized implicit schemes are proposed to approximate such the nonlinear equation
with continuous/discontinuous coefficients. The stabilities and convergences of the two schemes are proved under several suitable assumptions,
and numerical examples show that the convergence orders of these two schemes
 are $1$ in both time and space.
Secondly, a nonlinear all-at-once system is derived based on the nonlinear implicit scheme,
which may suitable for parallel computations.
Newton's method, whose initial value is obtained by interpolating the solution
of the linearized implicit scheme on the coarse space, is chosen to solve such the nonlinear all-at-once system.
To accelerate the speed of solving the Jacobian equations appeared in Newton's method,
a robust preconditioner is developed and analyzed. Numerical examples are reported to demonstrate
the effectiveness of our proposed preconditioner.
Meanwhile, they also imply that such the initial guess for Newton's method is more suitable.
\end{abstract}

\begin{keyword} Nonlinear tempered fractional diffusion equation\sep All-at-once system; Newton's method\sep
Krylov subspace method\sep Toeplitz matrix\sep banded Toeplitz preconditioner
\MSC[2010] 65L05\sep 65N22\sep 65F10
\end{keyword}

\end{frontmatter}


\section{Introduction}
\label{sec1}

In this work, we mainly focus on solving the all-at-once system
arising from the nonlinear tempered fractional diffusion equation (NL-TFDE):
\begin{equation}
\begin{cases}
\frac{\partial u(x,t)}{\partial t}
= d_{+}(x) \sideset{_a}{^{\alpha,\lambda}_x}{\mathop{\mathbf{D}}} u(x,t)
+ d_{-}(x) \sideset{_x}{^{\alpha,\lambda}_b}{\mathop{\mathbf{D}}} u(x,t)
+ f(u(x,t), x, t), & (x,t) \in [a, b] \times (0, T],\\ 
u(a,t) = u(b,t) = 0, & 0 \leq t \leq T, \\
u(x,0) = u_{0}(x), & a \leq x \leq b,
\end{cases}
\label{eq1.1}
\end{equation}
where $\alpha \in (1,2)$, $\lambda \geq 0$, 
$d_{+}(x) \geq d_{-}(x) > 0$
and $f(u(x,t), x, t)$ is a source term which satisfies the Lipschitz condition:
\begin{equation*}
\mid f(r_1, x, t) - f(r_2, x, t) \mid \leq L \mid r_1 - r_2 \mid, ~\textrm{for~all}~r_1, r_2~\textrm{over}~[a,b] \times [0,T].
\end{equation*}
The variants of the left and right Riemann-Liouville tempered fractional
derivatives are respectively defined as \cite{cartea2007fluid,baeumer2010tempered,meerschaert2012stochastic}
\begin{equation}
\sideset{_a}{^{\alpha,\lambda}_x}{\mathop{\mathbf{D}}} u(x,t)
= \sideset{_a}{^{\alpha,\lambda}_x}{\mathop{D}} u(x,t)
- \alpha \lambda^{\alpha - 1} \frac{\partial u(x,t)}{\partial x} - \lambda^\alpha u(x,t),
\label{eq1.2}
\end{equation}
and
\begin{equation}
\sideset{_x}{^{\alpha,\lambda}_b}{\mathop{\mathbf{D}}} u(x,t)
= \sideset{_x}{^{\alpha,\lambda}_b}{\mathop{D}} u(x,t)
+ \alpha \lambda^{\alpha - 1} \frac{\partial u(x,t)}{\partial x} - \lambda^\alpha u(x,t),
\label{eq1.3}
\end{equation}
where $\sideset{_a}{^{\alpha,\lambda}_x}{\mathop{D}} u(x,t)$ and
$\sideset{_x}{^{\alpha,\lambda}_b}{\mathop{D}} u(x,t)$ are the
left and right Riemann-Liouville tempered fractional derivatives
defined respectively as \cite{cartea2007fluid,baeumer2010tempered}
\begin{equation*}
\sideset{_a}{^{\alpha,\lambda}_x}{\mathop{D}} u(x,t)
= \frac{e^{-\lambda x}}{\Gamma(2 - \alpha)} \frac{\partial^2}{\partial x^2} \int_a^x
\frac{e^{\lambda \xi} u(\xi,t)}{(x - \xi)^{1 - \alpha}} d \xi,
\end{equation*}
\begin{equation*}
\sideset{_x}{^{\alpha,\lambda}_b}{\mathop{D}} u(x,t)
= \frac{e^{\lambda x}}{\Gamma(2 - \alpha)} \frac{\partial^2}{\partial x^2} \int_x^b
\frac{e^{-\lambda \xi} u(\xi,t)}{(\xi - x)^{1 - \alpha}} d \xi.
\end{equation*}
If $\lambda = 0$, they reduce to the Riemann-Liouville fractional derivatives \cite{podlubny1998fractional}.

Tempered fractional diffusion equations (TFDEs) are exponentially tempered extension
of fractional diffusion equations.
In recent several decades, the TFDEs are widely used across various fields,
such as statistical physics \cite{cartea2007fluid,chakrabarty2011tempered,zheng2015numerical},
finance \cite{carr2002fine,carr2003stochastic,wang2015circulant,zhang2016numerical} and geophysics \cite{baeumer2010tempered,meerschaert2008tempered,metzler2004restaurant,zhang2011gaussian,zhang2012linking}.
Unfortunately, it is difficult to obtain the analytical solutions of TFDEs,
or the obtained analytical solutions are less practical.
Hence, numerical methods such as finite difference method \cite{zhao2018fast,gu2017fast}
and finite element method \cite{li2018fast} become important approaches to solve TFDEs.
There are limited works addressing the finite difference schemes for the TFDEs.
Baeumera and Meerschaert \cite{baeumer2010tempered} provided finite difference
and particle tracking methods for solving the TFDEs on a bounded interval.
The stability and second-order accuracy of the resulted schemes are discussed.
Cartea and del-Castillo-Negrete \cite{cartea2007fractional} proposed a general
finite difference scheme to numerically solve a Black-Merton-Scholes model with tempered fractional derivatives.
Marom and Momoniat \cite{marom2009comparison} compared the numerical solutions
of three fractional partial differential equations (FDEs)
with tempered fractional derivatives that occur in finance.
However, the stabilities of their proposed schemes are not proved.
Recently, Li and Deng \cite{li2016high} derived a series of high order difference
approximations (called tempered-WSGD operators) for the tempered fractional calculus.
They also used such operators to numerically solve the TFDE,
and the stability and convergence of the obtained numerical schemes are proved.

Similar to the fractional derivatives, the tempered fractional derivatives are nonlocal.
Thus the discretized systems for TFDEs usually accompany a full (or dense) coefficient matrix.
Traditional methods (e.g., Gaussian elimination) to solve such systems need computational cost
is of $\mathcal{O}(N^3)$ and storage requirement is of $\mathcal{O}(N^2)$,
where $N$ is the number of space grid points.
Fortunately, the coefficient matrix always holds a Toeplitz-like structure.
It is well known that Toeplitz matrices possess great structures and properties,
and their matrix-vector multiplications can be computed in $\mathcal{O}(N\log N)$ operations
via fast Fourier transform (FFT) \cite{ng2004iterative,chan2007introduction}.
With this truth, the memory requirement and computational cost
of Krylov subspace methods are $\mathcal{O}(N)$ and $\mathcal{O}(N\log N)$, respectively.
However, the convergence rate of the Krylov subspace methods will be slow,
if the coefficient matrix is ill-conditioned.
To address this problem, Wang et al. \cite{wang2015circulant} proposed a circulant preconditioned
generalized minimal residual method (PGMRES) to solve the discretized linear system,
whose computational cost is of $\mathcal{O}(N\log N)$.
Lei et al. \cite{lei2018fast} proposed fast solution algorithms for solving TFDEs
in one-dimensional (1D) and two-dimensional (2D).
In their article, for 1D case, a circulant preconditioned iterative method and a fast-direct method are developed,
and the computational complexity of both methods are $\mathcal{O}(N\log N)$ in each time step.
For 2D case, such two methods were extended to fast solve their alternating direction implicit (ADI) scheme,
and the complexity of both methods are $\mathcal{O}(N^2 \log N)$ in each time step.
For many other studies about Toeplitz-like systems, 
see \cite{qu2016cscs,gu2015k,gu2015strang,zhao2012dct} and the references therein.

Actually, all the aforementioned fast implementations for TFDEs are developed based on the time-stepping schemes,
which are not suitable for parallel computations.
If all the time steps are stacked in a vector, the all-at-once system is obtained
and it is suitable for parallel computations, see \cite{Martin2017time, wu2018parareal}.
To the best of our knowledge, such the system arising from the FDEs
or the partial differential equations have been studied by many researchers
\cite{Martin2015time,Banjai2012multistep, McDonald2018pde, ke2015fast,lu2015fast,huang2017fast,lu2018approximate,zhao2018limited}.
However, the all-at-once system arising from the TFDEs is less studied.
In this work, a preconditioning technique is designed for such the system arising from the NL-TFDE \eqref{eq1.1}.
The rest of this paper is organized as follows: in Section \ref{sec2},
the nonlinear and linearized implicit schemes are derived by utilizing the finite difference method.
Then, the nonlinear all-at-once system is obtained from the nonlinear implicit scheme.
The stabilities and convergences of such two schemes are analyzed in Section \ref{sec3}.
A preconditioning technique is designed in Section \ref{sec4} to accelerate solving such the all-at-once system.
In Section \ref{sec5}, numerical examples are provided to illustrate the first-order convergences
of the two implicit schemes and show the performance of our preconditioning strategy
for solving such the system. Concluding remarks are given in Section \ref{sec6}.

\section{Two implicit schemes and all-at-once system}
\label{sec2}\quad
In this section, the nonlinear and linearized implicit schemes are proposed to approach Eq. \eqref{eq1.1}.
Then, the all-at-once system is obtained from the nonlinear one.

\subsection{Two implicit schemes}
\label{sec2.1}\quad
In order to derive the proposed schemes,
we first introduce the mesh $\bar{\omega}_{h \tau} = \bar{\omega}_{h} \times \bar{\omega}_{\tau}$,
where
$\bar{\omega}_{h} = \{ x_i = a + ih, ~i = 0, 1, \cdots, N; ~x_0 = a, x_N = b \}$ and
$\bar{\omega}_{\tau} = \{ t_j = j \tau, ~j = 0, 1, \cdots, M; ~t_M = T \}$.
Let $u_i^j$ represents the numerical approximation of $u(x_i,t_j)$.
Then the variants of the Riemann-Liouville tempered fractional derivatives defined in
Eqs. \eqref{eq1.2}-\eqref{eq1.3} for $(x,t) = (x_i,t_j)$ can be approximated respectively
as \cite{sabzikar2015tempered, li2016high}:
\begin{equation}
\sideset{_a}{^{\alpha,\lambda}_x}{\mathop{\mathbf{D}}} u(x,t) |_{(x,t) = (x_i, t_j)}
= \frac{1}{h^{\alpha}} \sum\limits_{k = 0}^{i + 1} g_k^{(\alpha)} u_{i - k + 1}^j
- \alpha \lambda^{\alpha - 1} \delta_x u_i^j + \mathcal{O}(h);
\label{eq2.1}
\end{equation}
\begin{equation}
\sideset{_x}{^{\alpha,\lambda}_b}{\mathop{\mathbf{D}}} u(x,t) |_{(x,t) = (x_i, t_j)}
= \frac{1}{h^{\alpha}} \sum\limits_{k = 0}^{N - i + 1} g_k^{(\alpha)} u_{i + k - 1}^j
+ \alpha \lambda^{\alpha - 1} \delta_x u_i^j + \mathcal{O}(h),
\label{eq2.2}
\end{equation}
where
\begin{equation*}
\delta_x u_i^j = \frac{u_i^j - u_{i - 1}^j}{h} \quad \mathrm{and} \quad
g_k^{(\alpha)} =
\begin{cases}
\tilde{g}_1^{(\alpha)} - e^{h \lambda} \left( 1 - e^{-h \lambda} \right)^{\alpha}, & k = 1,\\
\tilde{g}_k^{(\alpha)} e^{-(k - 1) h \lambda}, & k \neq 1
\end{cases}
\end{equation*}
with $\tilde{g}_k^{(\alpha)} = (-1)^k \binom{\alpha}{k}~(k \geq 0)$.

As for the time discretization, the backward Euler method is used.
Combining Eqs. \eqref{eq2.1} and \eqref{eq2.2},
the following first-order nonlinear implicit Euler scheme (NL-IES) is obtained:
\begin{equation}
\begin{split}
& u_i^j - w_1 \left( d_{+, i} \sum\limits_{k = 0}^{i + 1} g_k^{(\alpha)} u_{i - k + 1}^j
+ d_{-, i}\sum\limits_{k = 0}^{N - i + 1} g_k^{(\alpha)} u_{i + k - 1}^j \right)
+ w_2 \left( d_{+, i} - d_{-, i} \right) \left( u_i^j - u_{i - 1}^j \right) \\
& = u_i^{j - 1} + \tau f_{u, i}^j,
\end{split}
\label{eq2.3}
\end{equation}
in which $w_1 = \frac{\tau}{h^{\alpha}}$, $w_2 = \frac{\alpha \lambda^{\alpha - 1} \tau}{h}$,
$d_{\pm, i} = d_{\pm}(x_i)$ and $f_{u, i}^j = f(u(x_i, t_j), x_i, t_j)$.
Applying the formula $f(u(x_i, t_j), x_i, t_j) = f(u(x_i, t_{j - 1}), x_i, t_{j - 1}) + \mathcal{O}(\tau)$
to Eq. \eqref{eq2.3} and omitting the small term,
it gets the first-order linearized implicit Euler scheme (L-IES):
\begin{equation}
\begin{split}
& u_i^j - w_1 \left( d_{+, i} \sum\limits_{k = 0}^{i + 1} g_k^{(\alpha)} u_{i - k + 1}^j
+ d_{-, i}\sum\limits_{k = 0}^{N - i + 1} g_k^{(\alpha)} u_{i + k - 1}^j \right)
+ w_2 \left( d_{+, i} - d_{-, i} \right) \left( u_i^j - u_{i - 1}^j \right) \\
& = u_i^{j - 1} + \tau f_{u, i}^{j - 1}.
\end{split}
\label{eq2.4}
\end{equation}
The stabilities and first-order convergences of schemes \eqref{eq2.3}-\eqref{eq2.4}
will be discussed in Section \ref{sec3}.

\subsection{The all-at-once system}
\label{sec2.2}\quad
Several auxiliary notations are introduced before deriving the all-at-once system:
$I$ and $\bm{0}$ represent the identity and zero matrices of suitable orders, respectively.
\begin{equation*}
\bm{u}^j = \left[ u_1^j, u_2^j, \cdots, u_{N - 1}^j \right]^T, \qquad
\bm{f}_u^j = \left[ f_{u, 1}^j, f_{u, 2}^j, \cdots, f_{u, {N - 1}}^j \right]^T,
\end{equation*}
\begin{equation*}
D_{\pm} = \textrm{diag} (d_{\pm, 1}, d_{\pm, 2}, \cdots, d_{\pm, N - 1}), \qquad
B = \textrm{tridiag}(-1, 1, 0),
\end{equation*}
\begin{equation*}
\bm{u} = \begin{bmatrix}
\bm{u}^1 \\
\bm{u}^2 \\
\vdots \\
\bm{u}^M
\end{bmatrix}, ~
\bm{f}(\bm{u}) = \begin{bmatrix}
\bm{f}_u^1 \\
\bm{f}_u^2 \\
\vdots \\
\bm{f}_u^M
\end{bmatrix},~
\bm{v} = \begin{bmatrix}
\bm{u}^0 \\
\bm{0} \\
\vdots \\
\bm{0}
\end{bmatrix},~
G = \begin{bmatrix}
 g_1^{(\alpha)} & g_0^{(\alpha)} & 0 & \cdots & 0 \\
 g_2^{(\alpha)} & g_1^{(\alpha)} & g_0^{(\alpha)} & \ddots & 0 \\
 \vdots & \ddots & \ddots &  \ddots & \vdots \\
 g_{N-2}^{(\alpha)} & \cdots & \ddots & \ddots & g_0^{(\alpha)} \\
 g_{N-1}^{(\alpha)} & g_{N-2}^{(\alpha)} & \cdots &  g_2^{(\alpha)} & g_1^{(\alpha)}
\end{bmatrix}.
\end{equation*}

In this work, the all-at-once system is derived based on Eq. \eqref{eq2.3},
 which can be expressed as:
\begin{equation}
\mathcal{A} \bm{u} = \tau \bm{f}(\bm{u}) + \bm{v},
\label{eq2.5}
\end{equation}
in which $\mathcal{A} = \textrm{blktridiag} (-I, A, \bm{0})$ is a bi-diagonal block matrix with
$A = I - w_1 \left( D_{+} G + D_{-} G^T \right) + w_2 \left( D_{+} - D_{-} \right) B$.
Obviously, $A$ is a Toeplitz-like matrix and its storage requirement is of $\mathcal{O}(N)$.

For this nonlinear all-at-once system, we prefer to utilize Newton's method \cite{kelley2003solving}.
Such the method requires to solve the equation with Jacobian matrix at each iterative step,
and the computation of solving these equations consumes most of the method.
Before applying Newton's method to solve the system \eqref{eq2.5}, two essential problems need to be addressed:
\begin{itemize}
	\item[1.] {How to find a good enough initial value?}
	\item[2.] {How to solve the Jacobian equations efficiently?}
\end{itemize}

Here, a strategy is provided to address these two problems.
For the first problem, the initial value of Newton's method is constructed by interpolating the solution
of L-IES (\ref{eq2.4}) on the coarse mesh. Numerical experiences in Section 5 show that
it is a good enough initial value. For the second problem, the Jacobian matrix of (\ref{eq2.5})
is a bi-diagonal block matrix. More precisely, such the matrix is the sum of a diagonal block matrix
and a bi-diagonal block matrix, whose blocks are Toeplitz-like matrices.
Based on this special structure, the preconditioned Krylov subspace methods,
such as preconditioned biconjugate gradient stabilized (PBiCGSTAB) method \cite{van1992bicgstab},
are employed to solve the Jacobian equations appeared in Newton's method.
The details will be discussed in Section \ref{sec4}.

\section{Stabilities and convergences of (2.3) and (2.4)}
\label{sec3}\quad
In this section, the stabilities and convergences of NL-IES \eqref{eq2.3} and L-IES \eqref{eq2.4} are studied.
Let $U_i^j$ be the approximation solution of $u_i^j$ in Eqs. \eqref{eq2.3} or \eqref{eq2.4} and
$e_i^j = U_i^j - u_i^j ~(i = 1, \cdots, N - 1; ~j = 1, \cdots, M)$ be the error satisfying equation
\begin{equation*}
\begin{split}
& e_i^j - w_1 \left( d_{+, i} \sum\limits_{k = 0}^{i + 1} g_k^{(\alpha)} e_{i - k + 1}^j
+ d_{-, i}\sum\limits_{k = 0}^{N - i + 1} g_k^{(\alpha)} e_{i + k - 1}^j \right)
+ w_2 \left( d_{+, i} - d_{-, i} \right) \left( e_i^j - e_{i - 1}^j \right) \\
& = e_i^{j - 1} + \tau \left( f_{U, i}^j - f_{u, i}^j \right) ~(\textrm{for}~\textrm{NL-IES})
\end{split}
\end{equation*}
or
\begin{equation*}
\begin{split}
& e_i^j - w_1 \left( d_{+, i} \sum\limits_{k = 0}^{i + 1} g_k^{(\alpha)} e_{i - k + 1}^j
+ d_{-, i}\sum\limits_{k = 0}^{N - i + 1} g_k^{(\alpha)} e_{i + k - 1}^j \right)
+ w_2 \left( d_{+, i} - d_{-, i} \right) \left( e_i^j - e_{i - 1}^j \right) \\
& = e_i^{j - 1} + \tau \left( f_{U, i}^{j - 1} - f_{u, i}^{j - 1} \right)~(\textrm{for}~\textrm{L-IES}),
\end{split}
\end{equation*}
in which $f_{U, i}^{k} = f(U_i^k, x_i, t_k)$.
To prove the stabilities and convergences of \eqref{eq2.3} and \eqref{eq2.4},
the following results given in \cite{sabzikar2015tempered, zhuang2009numerical} are required.
\begin{lemma}(\cite{sabzikar2015tempered} )
The coefficients $g_k^{(\alpha)}$, for $k = 0, 1, \cdots$, satisfy:
\begin{equation*}
\begin{cases}
g_1^{(\alpha)} < 0,~ g_k^{(\alpha)} > 0~ (\textrm{for}~ k \neq 1); \\
\sum\limits_{k = 0}^{\infty} g_k^{(\alpha)} = 0,
~\sum\limits_{k = 0}^{j} g_k^{(\alpha)} < 0~ (\textrm{for}~ j \geq 1).
\end{cases}
\end{equation*}
\label{lemma3.1}
\end{lemma}
\begin{lemma}(\cite{zhuang2009numerical}, discrete Gronwall inequality)
Suppose that $\tilde{f}_k \geq 0,~\eta_k \geq 0~(k = 0, 1, \cdots)$, and
\begin{equation*}
\eta_{k + 1} \leq \rho \eta_k + \tau \tilde{f_k},~\rho = 1 + C_0 \tau,  ~\eta_0 = 0,
\end{equation*}
where $C_0 \geq 0$ is a constant, then
\begin{equation*}
\eta_{k + 1} \leq \exp(C_0 t_k) \sum\limits_{j = 0}^{k} \tau \tilde{f}_j.
\end{equation*}
\label{lemma3.2}
\end{lemma}

\subsection{The stabilities of (2.3) and (2.4)}
\label{sec3.1}\quad
Denote $E^j = \left[ e_1^j, e_2^j, \cdots, e_{N - 1}^j \right]^T$,
and assume that
\begin{equation*}
\left\| E^j \right\|_\infty = \mid e_{\ell_j} \mid
= \max\limits_{1 \leq \ell \leq N - 1} \mid e_{\ell}^j \mid~(0 \leq j \leq M, 1 \leq \ell_j \leq N - 1).
\end{equation*}
Then the next theorem is established.
\begin{theorem}
Suppose $d_{+}(x) \geq d_{-}(x) \geq 0$, then the L-IES \eqref{eq2.4} is stable, and we have
\begin{equation*}
\left\| E^j \right\|_\infty \leq \exp(TL) \left\| E^0 \right\|_\infty, ~\textrm{for}~ j= 1, \cdots, M.
\end{equation*}
\label{th3.1}
\end{theorem}
\textbf{Proof.}
From Lemma \ref{lemma3.1},
$\sum\limits_{k = 0}^{\ell_j + 1} g_k^{(\alpha)} < 0$ and
$\sum\limits_{k = 0}^{N - \ell_j + 1} g_k^{(\alpha)} < 0$. Then
\begin{equation*}
\begin{split}
\mid e_{\ell_j}^j \mid
& \leq \left[ 1 - w_1 \left( d_{+, \ell_j} \sum\limits_{k = 0}^{\ell_j + 1} g_k^{(\alpha)}
+ d_{-, \ell_j} \sum\limits_{k = 0}^{N - \ell_j + 1} g_k^{(\alpha)} \right) \right] \mid e_{\ell_j}^j \mid
 + w_2 \left( d_{+, \ell_j} - d_{-, \ell_j} \right) \times \\
& \quad \left( \mid e_{\ell_j}^j \mid - \mid e_{\ell_j}^j \mid \right) \\
& \leq \mid e_{\ell_j}^j \mid - w_1 \left( d_{+, \ell_j} \sum\limits_{k = 0}^{\ell_j + 1} g_k^{(\alpha)}
\mid e_{\ell_j - k + 1}^j \mid
+ d_{-, \ell_j} \sum\limits_{k = 0}^{N - \ell_j + 1} g_k^{(\alpha)} \mid e_{\ell_j + k - 1}^j \mid \right) \\
&\quad + w_2 \left( d_{+, \ell_j} - d_{-, \ell_j} \right) \left( \mid e_{\ell_j}^j \mid - \mid e_{\ell_j - 1}^j \mid \right) \\
& = \mid e_{\ell_j}^j \mid + \mid -w_1 d_{+, \ell_j} g_1^{(\alpha)} e_{\ell_j}^j \mid
+ \mid -w_1 d_{-, \ell_j} g_1^{(\alpha)} e_{\ell_j}^j \mid
+ \mid w_2 \left( d_{+, \ell_j} - d_{-, \ell_j} \right) e_{\ell_j}^j \mid \\
&\quad - \sum\limits_{k = 0, k \neq 1}^{\ell_j + 1} \mid -w_1 d_{+, \ell_j} g_k^{(\alpha)} e_{\ell_j - k + 1}^j \mid
- \sum\limits_{k = 0, k \neq 1}^{N - \ell_j + 1} \mid -w_1 d_{-, \ell_j} g_k^{(\alpha)} e_{\ell_j + k - 1}^j \mid \\
&\quad - \mid -w_2 \left( d_{+, \ell_j} - d_{-, \ell_j} \right) e_{\ell_j - 1}^j \mid \\
& = \mid e_{\ell_j}^j - w_1 \left( d_{+, \ell_j} g_1^{(\alpha)} e_{\ell_j}^j
+ d_{-, \ell_j} g_1^{(\alpha)} e_{\ell_j}^j \right) + w_2 \left( d_{+, \ell_j} - d_{-, \ell_j} \right) e_{\ell_j}^j \mid \\
&\quad - \sum\limits_{k = 0, k \neq 1}^{\ell_j + 1} \mid -w_1 d_{+, \ell_j} g_k^{(\alpha)} e_{\ell_j - k + 1}^j \mid
- \sum\limits_{k = 0, k \neq 1}^{N - \ell_j + 1} \mid -w_1 d_{-, \ell_j} g_k^{(\alpha)} e_{\ell_j + k - 1}^j \mid \\
&\quad - \mid -w_2 \left( d_{+, \ell_j} - d_{-, \ell_j} \right) e_{\ell_j - 1}^j \mid \\
& \leq \mid e_{\ell_j}^j - w_1 \left( d_{+, \ell_j} \sum\limits_{k = 0}^{\ell_j + 1}
g_k^{(\alpha)} e_{\ell_j - k + 1}^j
+ d_{-, \ell_j} \sum\limits_{k = 0}^{N - \ell_j + 1} g_k^{(\alpha)} e_{\ell_j + k - 1}^j \right) \\
&\quad + w_2 \left( d_{+, \ell_j} - d_{-, \ell_j} \right) \left( e_{\ell_j}^j - e_{\ell_j - 1}^j \right) \mid \\
& = \mid e_{\ell_j}^{j - 1} + \tau \left( f_{U, \ell_j}^{j - 1} - f_{u, \ell_j}^{j - 1} \right) \mid
\leq \left(1 + \tau L \right) \mid e_{\ell_j}^{j - 1} \mid \\
& \leq \left(1 + \tau L \right) \left\| E^{j - 1} \right\|_\infty.
\end{split}
\end{equation*}

The above inequality implies $\left\| E^{j} \right\|_\infty \leq \left(1 + \tau L \right) \left\| E^{j - 1} \right\|_\infty$.
Then
\begin{equation*}
\left\| E^{j} \right\|_\infty \leq \left(1 + \tau L \right)^j \left\| E^{0} \right\|_\infty
\leq \exp(TL) \left\| E^{0} \right\|_\infty.
\end{equation*} \hfill$\Box$

From the above proof, it can be find that if $\tau L < 1$, the following result is true.
\begin{corollary}
Suppose $d_{+}(x) \geq d_{-}(x) \geq 0$ and $\tau L < 1$, then the NL-IES \eqref{eq2.3} is stable,
and it obtains
\begin{equation*}
\left\| E^k \right\|_\infty \leq C_1 \left\| E^0 \right\|_\infty, ~\textrm{for}~k= 1, \cdots, M,
\end{equation*}
where $C_1$ is a positive constant.
\label{coro1}
\end{corollary}
\textbf{Proof.}
Based on the proof of Theorem \ref{th3.1}, it yields
\begin{equation*}
\left(1 - \tau L \right) \left\| E^{j} \right\|_\infty \leq \left\| E^{j - 1} \right\|_\infty.
\end{equation*}
Summing up for $j$ from $1$ to $k$ and using Lemma 3.4 in \cite{zhao2018fast}, it gets
\begin{equation*}
\left\| E^{k} \right\|_\infty \leq \frac{\exp(\frac{TL}{1 - \tau L})}{1 - \tau L}\left\| E^{0} \right\|_\infty.
\end{equation*}
Note that
\begin{equation*}
\lim\limits_{\tau \to 0} \frac{\exp(\frac{TL}{1 - \tau L})}{1 - \tau L}
= \exp(TL).
\end{equation*}
Hence, there is a positive constant $C_1$ such that
\begin{equation*}
\frac{\exp(\frac{TL}{1 - \tau L})}{1 - \tau L} \leq C_1,
\end{equation*}
thereby $\left\| E^k \right\|_\infty \leq C_1 \left\| E^0 \right\|_\infty, ~\textrm{for}~ k= 1, \cdots, M$.
\hfill$\Box$

In Corollary \ref{coro1}, it is worth to notice that the assumption $\tau < \frac{1}{L}$
is independent of the spatial size $h$.
Actually, the smaller time step size $\tau$ is, the easier such assumption can be satisfied.
\subsection{The convergences of (2.3) and (2.4)}
\label{sec3.2}\quad
In this subsection, the convergences of \eqref{eq2.3} and \eqref{eq2.4} are studied.
Let $\xi_i^j = u(x_i, t_j) - u_i^j$ satisfies
\begin{equation*}
\begin{split}
& \xi_i^j - w_1 \left( d_{+, i} \sum\limits_{k = 0}^{i + 1} g_k^{(\alpha)} \xi_{i - k + 1}^j
+ d_{-, i}\sum\limits_{k = 0}^{N - i + 1} g_k^{(\alpha)} \xi_{i + k - 1}^j \right)
+ w_2 \left( d_{+, i} - d_{-, i} \right) \left( \xi_i^j - \xi_{i - 1}^j \right) \\
& = \xi_i^{j - 1} + \tau \left( f(u(x_i, t_j), x_i, t_j) - f_{u, i}^j \right) + R_i^j
~(\textrm{for}~\textrm{NL-IES})
\end{split}
\end{equation*}
or
\begin{equation*}
\begin{split}
& \xi_i^j - w_1 \left( d_{+, i} \sum\limits_{k = 0}^{i + 1} g_k^{(\alpha)} \xi_{i - k + 1}^j
+ d_{-, i}\sum\limits_{k = 0}^{N - i + 1} g_k^{(\alpha)} \xi_{i + k - 1}^j \right)
+ w_2 \left( d_{+, i} - d_{-, i} \right) \left( \xi_i^j - \xi_{i - 1}^j \right) \\
& = \xi_i^{j - 1} + \tau \left(  f(u(x_i, t_{j - 1}), x_i, t_{j - 1}) - f_{u, i}^{j - 1} \right) + R_i^j ~(\textrm{for}~\textrm{L-IES}),
\end{split}
\end{equation*}
where $\left| R_i^j \right| \leq C_2 \left( \tau^2 + \tau h \right)$ ($C_2$ is a positive constant).
Denote $\bm{\xi^j} = \left[ \xi_1^j, \xi_2^j, \cdots, \xi_{N - 1}^j \right]^T$ and
$\left\| \bm{\xi}^j \right\|_\infty = \mid \xi_{\ell_j} \mid
= \max\limits_{1 \leq \ell \leq N - 1} \mid \xi_{\ell}^j \mid~(0 \leq j \leq M, 1 \leq \ell_j \leq N - 1)$.
Similar to the proof of Theorem \ref{th3.1}, the following theorem
about the convergence of  L-IES can be established.
\begin{theorem}
Assume that $d_{+}(x) \geq d_{-}(x) \geq 0$ and the problem \eqref{eq1.1}
has a sufficiently smooth solution $u(x,t)$.
$u_i^j$ is the numerical solution of \eqref{eq2.4}.
Then there is a positive constant $C$ such that
\begin{equation*}
\left\| \bm{\xi}^j \right\|_\infty \leq C (\tau + h), ~j = 1, 2, \cdots, M.
\end{equation*}
\label{th3.2}
\end{theorem}
\textbf{Proof.}
Same technique in Theorem \ref{th3.1} is utilized, then it yields
\begin{equation*}
\begin{split}
\left\| \bm{\xi}^{j} \right\|_\infty = \mid \xi_{\ell_j}^j \mid
& \leq
\mid \xi_{\ell_j}^{j - 1} + \tau \left( f(u(x_{\ell_j}, t_{j - 1}), x_{\ell_j}, t_{j - 1}) - f_{u, \ell_j}^{j - 1} \right)
+ R_{\ell_j}^j \mid \\
& \leq \left( 1 + \tau L \right) \mid \xi_{\ell_j}^{j - 1} \mid + C_2 \left( \tau^2 + \tau h \right) \\
& \leq \left( 1 + \tau L \right) \left\| \bm{\xi}^{j - 1} \right\|_\infty + C_2 \left( \tau^2 + \tau h \right) .
\end{split}
\end{equation*}
Using Lemma \ref{lemma3.2}, it gets
\begin{equation*}
\left\| \bm{\xi}^{j} \right\|_\infty
\leq \exp(TL) T C_2 \left( \tau + h \right) \leq C \left( \tau + h \right).
\end{equation*}
\hfill$\Box$

\begin{corollary}
Assume that $d_{+}(x) \geq d_{-}(x) \geq 0$, $\tau L < 1$
and the problem \eqref{eq1.1} has a sufficiently smooth solution $u(x,t)$.
$u_i^j$ is the numerical solution of \eqref{eq2.3}.
Then
\begin{equation*}
\left\| \bm{\xi}^j \right\|_\infty \leq C (\tau + h), ~j = 1, 2, \cdots, M,
\end{equation*}
where $C$ is a positive constant.
\label{coro2}
\end{corollary}
\textbf{Proof.}
According to Corollary \ref{coro1}, we have
\begin{equation*}
\left(1 - \tau L \right) \left\| \bm{\xi}^{j} \right\|_\infty \leq \left\| \bm{\xi}^{j - 1} \right\|_\infty
+ C_2 \left( \tau^2 + \tau h \right).
\end{equation*}
Similarly, it arrives
\begin{equation*}
\left\| \bm{\xi}^{j} \right\|_\infty \leq \frac{\exp(\frac{TL}{1 - \tau L})}{1 - \tau L}
T C_2 \left( \tau + h \right)
\leq C \left( \tau + h \right), ~j = 1, 2, \cdots, M,
\end{equation*}
in which $\frac{\exp(\frac{TL}{1 - \tau L})}{1 - \tau L} \leq C_1$ is employed.
\hfill$\Box$

It is interesting to note that if $\xi_i^j$ represents the error between NL-IES \eqref{eq2.3} and L-IES \eqref{eq2.4},
then it also satisfies  $\left\| \bm{\xi}^j \right\|_\infty \leq C (\tau + h)~(j = 1, 2, \cdots, M)$.
This can be proved easily through the condition $f_{u,i}^{j} = f_{u, i}^{j - 1} + \mathcal{O}(\tau)$.
In the next section, fast implementations are designed to solve \eqref{eq2.5}.

\section{The preconditioned iterative method}
\label{sec4}\quad
The Jacobian matrix of \eqref{eq2.5} can be treated as the sum of a diagonal block matrix
and a block bi-diagonal matrix with Toeplitz-like blocks,
then its matrix-vector multiplication can be done by FFT in $\mathcal{O}(MN \log MN)$ operations.
Such the technique truly reduces the computational cost of a Krylov subspace method,
such as the biconjugate gradient stabilized (BiCGSTAB) method,
but the convergence rate of this method is slow when the Jacobi matrix is very ill-conditioned.
In order to speed up the convergence rate of the Krylov subspace method,
a preconditioner $P_{\ell} = \textrm{blktridiag} (-I, A_{\ell}, \bm{0})~(\ell > 2)$
is proposed and analyzed in this section, in which
\begin{equation*}
A_{\ell} = I - w_1 \left( D_{+} G_{\ell} + D_{-} G_{\ell}^T \right) + w_2 \left( D_{+} - D_{-} \right) B ~
\textrm{with}~
G_{\ell} = \begin{bmatrix}
	g_1^{(\alpha)} & g_0^{(\alpha)} & & & & \\
	g_2^{(\alpha)} & g_1^{(\alpha)} & g_0^{(\alpha)} & & & \\
	\vdots & \ddots & \ddots & \ddots & & \\
	g_{\ell}^{(\alpha)} & \cdots & g_2^{(\alpha)} & g_1^{(\alpha)} & g_0^{(\alpha)} & \\
	 & \ddots & \ddots & \ddots & \ddots & \ddots \\
	 & & g_{\ell}^{(\alpha)} & \cdots & g_2^{(\alpha)} & g_1^{(\alpha)}
\end{bmatrix}.
\end{equation*}

Noticing the properties of $g_{k}^{(\alpha)}$ given in Lemma \ref{lemma3.1},
the following result about $P_{\ell}$ is true.
\begin{theorem}
The preconditioner $P_{\ell}$ is a nonsingular matrix.
\label{th4.1}
\end{theorem}
\textbf{Proof.}
Obviously, we only need to proof the nonsingularity of $A_{\ell}$. Let
\begin{equation*}
H = \frac{A_{\ell} + A_{\ell}^T}{2}
= I - \frac{\omega_1}{2} \left( D_{+} + D_{-} \right) \left( G_{\ell} + G_{\ell}^T \right)
+ \frac{\omega_2}{2} \left( D_{+} - D_{-} \right) \left( B + B^T \right)
\end{equation*}
and $\theta$ be the arbitrary eigenvalue of $H$.
According to Lemma \ref{lemma3.1} and Gershgorin circle theorem \cite{varga2010gervsgorin},
it arrives at
\begin{equation*}
\begin{split}
& \left| \theta - \left[ 1 - \omega_1 \left(d_{+,i} + d_{-,i}\right) g_1^{(\alpha)}
+ \omega_2 \left(d_{+,i} - d_{-,i}\right) \right] \right| \\
& \leq r_i =
\left| -\omega_1 \left(d_{+,i} + d_{-,i}\right) \left(g_0^{(\alpha)} + g_2^{(\alpha)}\right)
- \omega_2 \left(d_{+,i} - d_{-,i}\right) \right|
+ \sum_{k = 3}^{\ell} \left|-\omega_1 \left(d_{+,i} + d_{-,i}\right) g_k^{(\alpha)}\right|\\
& \leq
\omega_1 \left(d_{+,i} + d_{-,i}\right) \sum_{k = 0,~k \neq 1}^{\ell} g_k^{(\alpha)}
+ \omega_2 \left(d_{+,i} - d_{-,i}\right)
< - \omega_1 \left(d_{+,i} + d_{-,i}\right) g_1^{(\alpha)}
+ \omega_2 \left(d_{+,i} - d_{-,i}\right).
\end{split}
\end{equation*}

This implies that all eigenvalues of $H$ are larger than 1. Then, the desired result is achieved.
\hfill$\Box$

Unfortunately, it is difficult to theoretically investigate the eigenvalue distributions
of the preconditioned Jacobian matrix, but we still can give some figures
to illustrate the clustering eigenvalue distributions of several specified preconditioned matrices
in the next section.
For convenience, let $\bm{u}^{(k + 1)}$ be the approximation of $\bm{u}$
obtained in the $k$-th Newton iterative step,
the Jacobian matrix in the $k$-th Newton iterative step is denoted as $J^k$.
With these auxiliary notations, the preconditioned Newton's method can be summarized as below:
\begin{algorithm}[H]
	\caption{Solve $\bm{u}$ from Eq. \eqref{eq2.5}}
	\begin{algorithmic}[1]
		\STATE {Given maximum iterative step $maxit$, tolerance $tol_{out}$
			and initial vector $\bm{u}^{(0)}$, which is obtained by interpolating the solution
			of L-IES \eqref{eq2.4} on the coarse grid (here $M = N = 16$)}
		\FOR {$k = 1, \cdots, maxit$}
		\STATE {Solve $J^k \bm{z} = -\bm{f}(\bm{u}^{(k)})$ via PBiCGSTAB method
			with preconditioner $P_{\ell}$ ($\ell = 8$ is chosen experimentally to
			balance the number of iterations and CPU time)}
		\STATE {$\bm{u}^{(k + 1)} = \bm{u}^{(k)} + \bm{z}$}
		\IF {$\|\bm{z}\|_2 \leq tol_{out}$}
		\STATE $\bm{u} = \bm{u}^{(k + 1)}$
		\STATE \textbf{break}
		\ENDIF
		\ENDFOR
	\end{algorithmic}
	\label{alg1}
\end{algorithm}

In fact, this algorithm can be viewed as a simple two-grid method,
our readers can refer to \cite{xu1994twogrid,xu1996twogrid,kim2018twogrid} for details.

\section{Numerical examples}
\label{sec5}\quad
Two numerical experiments presented in this section have a two-fold objective.
On the one hand, they illustrate that the convergence orders
of our two implicit schemes \eqref{eq2.3}-\eqref{eq2.4} are 1.
On the other hand, they show the performance of the preconditioner $P_{\ell}$ proposed in Section \ref{sec4}.
In Algorithm \ref{alg1}, for generating the initial guess $\bm{u}^{(0)}$,
the MATLAB build-in function ``$\mathtt{interp2}$" is employed in this work.
The $maxit$ and $tol_{out}$ in Algorithm \ref{alg1} are fixed as $100$ and $10^{-12}$, respectively.
For the PBiCGSTAB method (or the BiCGSTAB method), it terminates if the relative residual error satisfies
$\frac{\| r^{k} \|_{2}}{\| r^{0} \|_{2}} \leq 10^{-6}$ or the iteration number is more than $1000$,
where $r^k$ is the residual vector of the linear system after $k$ iteration,
and the initial guess of the PBiCGSTAB method (or the BiCGSTAB method) is chosen as the zero vector.
''$\mathcal{P}$" represents that our proposed preconditioned iterative method
in Section \ref{sec4} is utilized to solve \eqref{eq2.5}.
``BS" (or ``$\mathcal{I}$") means that the PBiCGSTAB method in Step 3 of Algorithm \ref{alg1}
is replaced by the MATLAB's backslash method (or the BiCGSTAB method).
Some other notations, which will appear in later, are given:
\begin{equation*}
Err(\tau,h) = \max_{0 \leq j \leq M} \| \bm{\xi}^{j} \|_{\infty},
\end{equation*}
\begin{equation*}
Order1 = \log_{2} Err(2\tau,h)/Err(\tau,h), \quad
Order2 = \log_{2} Err(\tau,2h)/Err(\tau,h).
\end{equation*}
``$\mathrm{Iter1}$" represents the number of iterations required by Algorithm \ref{alg1}.
``$\mathrm{Iter2}$" denotes the average number of iterations required
by the PBiCGSTAB method (or the BiCGSTAB method) in Algorithm \ref{alg1},
i.e.,
\begin{equation*}
\mathrm{Iter2} = \sum\limits_{m = 1}^{\mathrm{Iter1}} \mathrm{Iter2} (m)/\mathrm{Iter1},
\end{equation*}
where $\mathrm{Iter2} (m)$ is the number of iterations required by such the method
in the $m$-th iterative step of Algorithm \ref{alg1}.
``$\textrm{Time}$" denotes the total CPU time in seconds for solving the system \eqref{eq2.5}.
``$\ddag$" means the maximum iterative step is reached but not convergence,
and``$\dag$" means out of memory.

All experiments were performed on a Windows 10 (64 bit) PC-Intel(R) Core(TM) i7-8700k CPU 3.70 GHz,
16 GB of RAM using MATLAB R2016a.

\begin{figure}[t]
	\centering
	\subfigure[Eigenvalues of $J^{0}$]
	{\includegraphics[width=3.0in,height=1.4in]{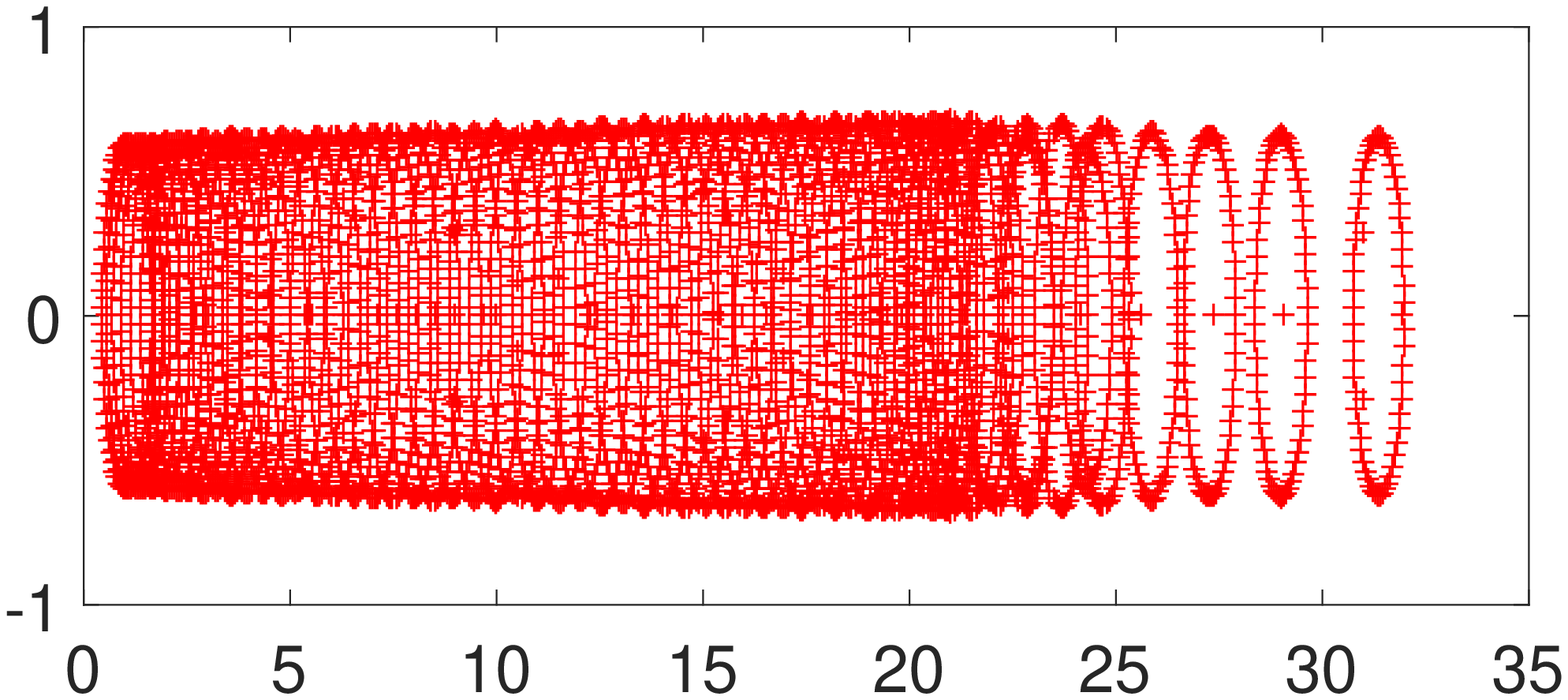}}\hspace{0.2in}
	\subfigure[Eigenvalues of $P_{\ell}^{-1} J^{0}$]
	{\includegraphics[width=3.0in,height=1.4in]{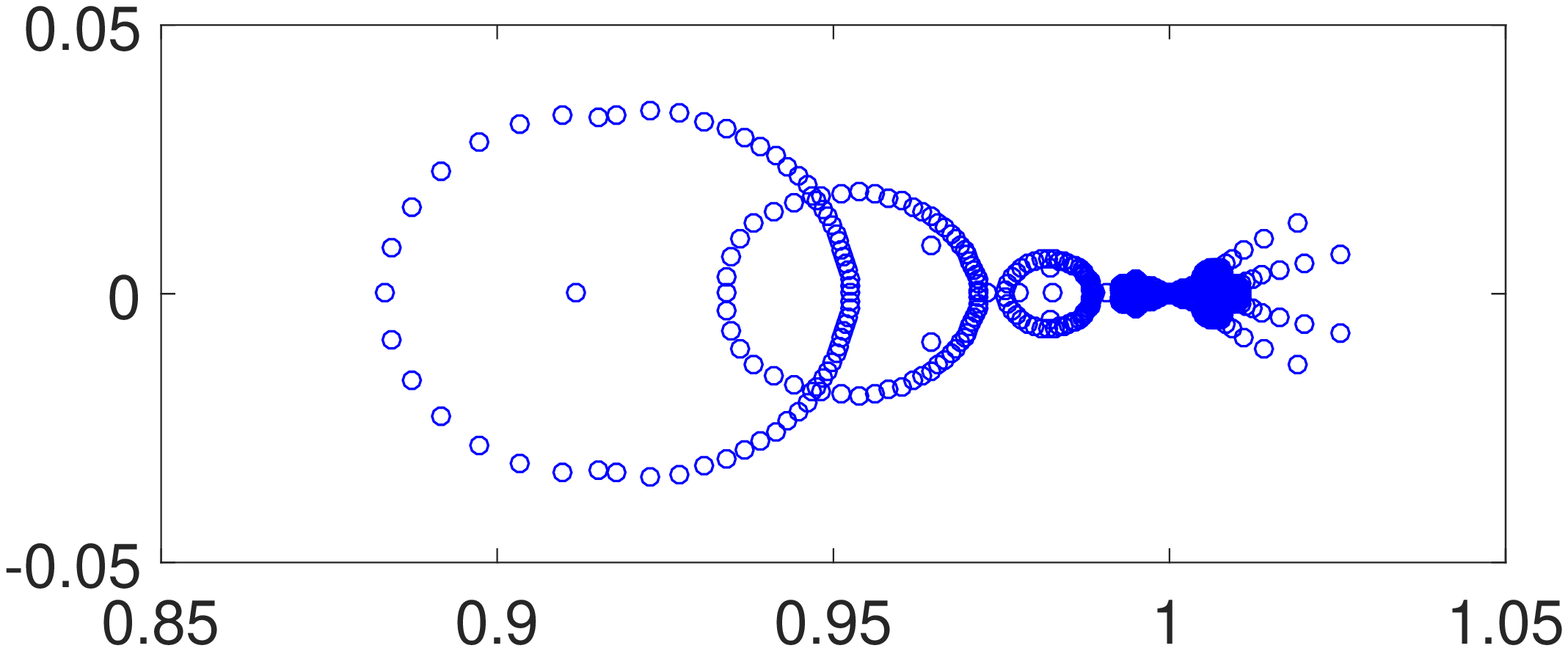}}\vspace{7.5em}
	\subfigure[Eigenvalues of $J^{0}$]
	{\includegraphics[width=3.0in,height=1.4in]{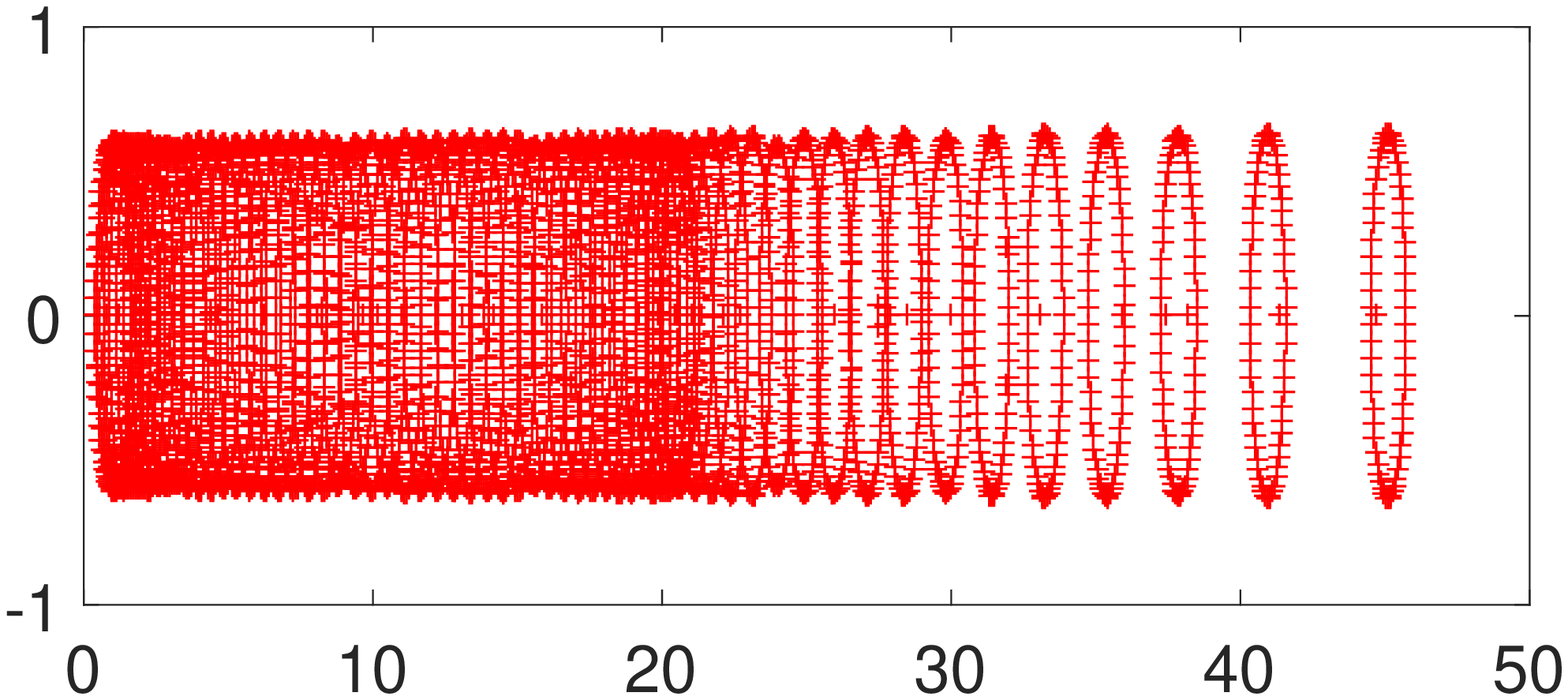}}\hspace{0.2in}
	\subfigure[Eigenvalues of $P_{\ell}^{-1} J^{0}$]
	{\includegraphics[width=3.0in,height=1.4in]{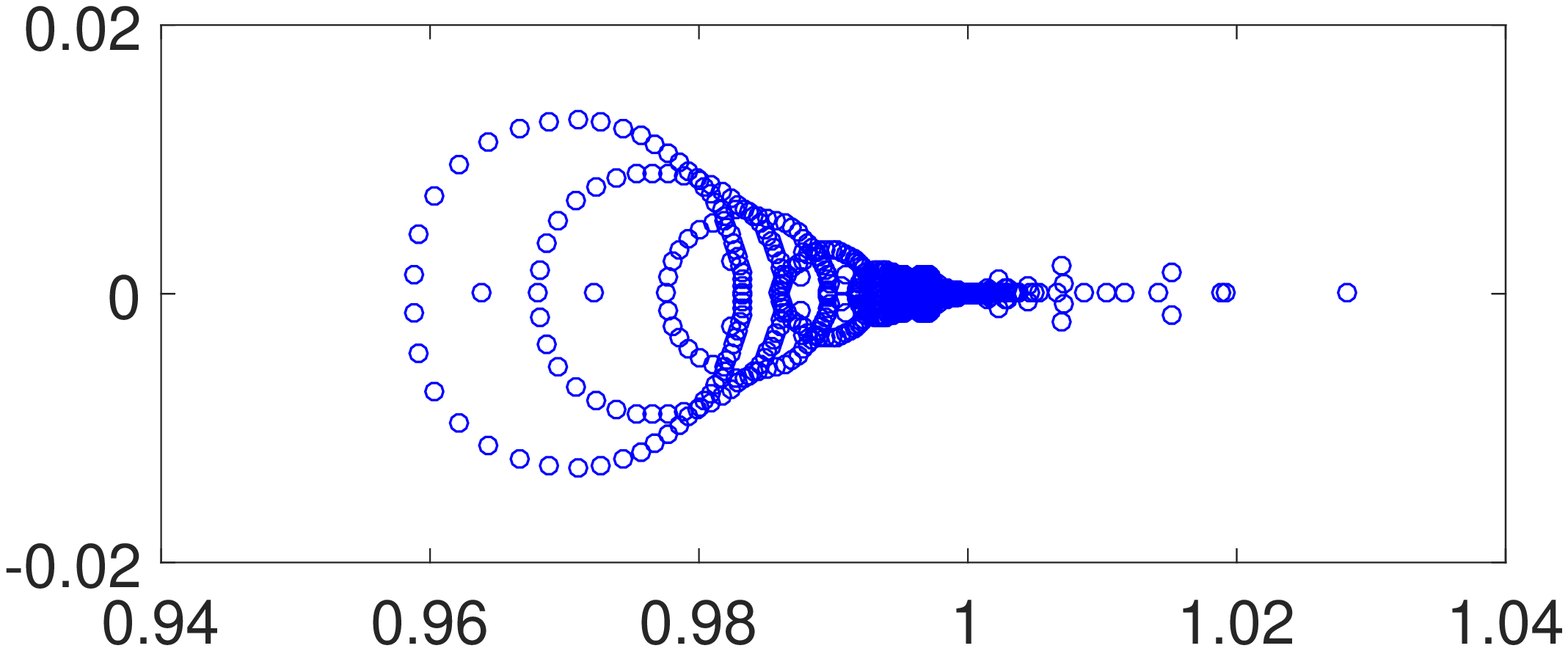}}
	\caption{Spectra of $J^{0}$ and $P_{\ell}^{-1} J^{0}$, when $\alpha = 1.5$, $M = N = 65$ in Example 1.
		Top row: $\lambda = 0$; Bottom row: $\lambda = 5$.}
	\label{fig2}
\end{figure}

\noindent \textbf{Example 1.} We consider Eq. \eqref{eq1.1} with $T = 1$, the initial value
$u(x,0) =  \left(\cos\pi x - 1 \right) \sin\pi x~(x \in [-1,1])$,
the nonlinear source term $f(u(x,t),x,t) = u(x,t) - 3 u(x,t)$
and the continuous coefficients $d_{+}(x) = 1.5 \exp(-x)$ and $d_{-}(x) = \exp(x)$.
Obviously, it is hard to obtain the exact solution of Eq. \eqref{eq1.1}.
Thus, the numerical solution computed from the finer mesh ($M = N = 1024$) is treated as the exact solution.
\begin{landscape}
\begin{table}[t]\scriptsize\tabcolsep=2.0pt
	\begin{center}
		\caption{The maximum norm errors and convergence orders for Example 1 with $h = 2^{-10}$.}
		\centering
		\begin{tabular}{cccccccccccccccccc}
			\hline
			& & \multicolumn{4}{c}{$\lambda = 0$} & \multicolumn{4}{c}{$\lambda = 1$}
			& \multicolumn{4}{c}{$\lambda = 5$} & \multicolumn{4}{c}{$\lambda = 10$} \\
			[-2pt] \cmidrule(lr){3-6} \cmidrule(lr){7-10} \cmidrule(lr){11-14} \cmidrule(lr){15-18}\\ [-11pt]
			& & \multicolumn{2}{c}{\rm{L-IES \eqref{eq2.4}}} & \multicolumn{2}{c}{\rm{NL-IES \eqref{eq2.3}}}
			& \multicolumn{2}{c}{\rm{L-IES \eqref{eq2.4}}} & \multicolumn{2}{c}{\rm{NL-IES \eqref{eq2.3}}}
			& \multicolumn{2}{c}{\rm{L-IES \eqref{eq2.4}}} & \multicolumn{2}{c}{\rm{NL-IES \eqref{eq2.3}}}
			& \multicolumn{2}{c}{\rm{L-IES \eqref{eq2.4}}} & \multicolumn{2}{c}{\rm{NL-IES \eqref{eq2.3}}} \\
			[-2pt] \cmidrule(lr){3-4} \cmidrule(lr){5-6} \cmidrule(lr){7-8} \cmidrule(lr){9-10}
			\cmidrule(lr){11-12} \cmidrule(lr){13-14} \cmidrule(lr){15-16} \cmidrule(lr){17-18} \\ [-11pt]
			$\alpha$ & $M$ & $Err(\tau,h)$ & $Order1$ & $Err(\tau,h)$ & $Order1$ & $Err(\tau,h)$ & $Order1$
			& $Err(\tau,h)$ & $Order1$ & $Err(\tau,h)$ & $Order1$ & $Err(\tau,h)$ & $Order1$
			& $Err(\tau,h)$ & $Order1$ & $Err(\tau,h)$ & $Order1$\\
\hline
1.1 & 64 & 8.0277E-02 & -- & 8.1498E-02 & -- & 1.1908E-02 & -- & 1.1390E-02 & -- & 6.6057E-03 & --
              & 6.3807E-03 & -- & 4.8182E-03 & -- & 4.6460E-03 & -- \\
    & 128 & 4.0105E-02 & 1.0012 & 4.1591E-02 & 0.9705 & 6.1196E-03 & 0.9604 & 5.8931E-03 & 0.9507 & 3.6702E-03 & 0.8479
              & 3.5897E-03 & 0.8299 & 2.7628E-03 & 0.8024 & 2.7029E-03 & 0.7815 \\
    & 256 & 1.7837E-02 & 1.1689 & 1.8690E-02 & 1.1540 & 2.8629E-03 & 1.0960 & 2.7890E-03 & 1.0793 & 1.8010E-03 & 1.0271
              & 1.7783E-03 & 1.0134 & 1.4017E-03 & 0.9790 & 1.3849E-03 & 0.9647 \\
    & 512 & 6.0553E-03 & 1.5586 & 6.3917E-03 & 1.5480 & 1.0254E-03 & 1.4813 & 1.0039E-03 & 1.4741 & 6.7680E-04 & 1.4120
              & 6.6838E-04 & 1.4118 & 5.2752E-04 & 1.4099 & 5.2131E-04 & 1.4096 \\
\hline
1.5 & 64 & 4.7666E-02 & -- & 5.3762E-02 & -- & 3.5875E-02 & -- & 4.3134E-02 & -- & 2.4027E-02 & --
              & 3.1109E-02 & -- & 1.8286E-02 & -- & 2.5143E-02 & -- \\
    & 128 & 3.6430E-02 & 0.3878 & 3.5974E-02 & 0.5796 & 1.9300E-02 & 0.8944 & 2.1332E-02 & 1.0158 & 1.1881E-02 & 1.0160
              & 1.5200E-02 & 1.0333 & 8.9830E-03 & 1.0255 & 1.2218E-02 & 1.0411 \\
    & 256 & 1.9950E-02 & 0.8687 & 1.9854E-02 & 0.8575 & 1.2203E-02 & 0.6614 & 1.2140E-02 & 0.8133 & 6.9339E-03 & 0.7769
              & 6.9000E-03 & 1.1394 & 4.7950E-03 & 0.9057 & 5.3432E-03 & 1.1932 \\
    & 512 & 7.5241E-03 & 1.4068 & 7.5004E-03 & 1.4044 & 5.2003E-03 & 1.2306 & 5.1850E-03 & 1.2274 & 3.2837E-03 & 1.0783
              & 3.2750E-03 & 1.0751 & 2.3859E-03 & 1.0070 & 2.3798E-03 & 1.1669 \\
\hline
1.9 & 64 & 1.1241E-01 & -- & 1.1930E-01 & -- & 1.1075E-01 & -- & 1.1761E-01 & -- & 1.0545E-01 & --
              & 1.1229E-01 & -- & 1.0186E-01 & -- & 1.0871E-01 & -- \\
    & 128 & 6.0869E-02 & 0.8850 & 6.3952E-02 & 0.8995 & 5.9987E-02 & 0.8846 & 6.3072E-02 & 0.8989 & 5.7145E-02 & 0.8839
              & 6.0268E-02 & 0.8978 & 5.5212E-02 & 0.8835 & 5.8362E-02 & 0.8974 \\
    & 256 & 2.8789E-02 & 1.0802 & 3.0258E-02 & 1.0797 & 2.8379E-02 & 1.0798 & 2.9832E-02 & 1.0801 & 2.6822E-02 & 1.0912
              & 2.8264E-02 & 1.0924 & 2.5776E-02 & 1.0990 & 2.7172E-02 & 1.1029 \\
    & 512 & 1.0164E-02 & 1.5020 & 1.0654E-02 & 1.5059 & 1.0020E-02 & 1.5019 & 1.0506E-02 & 1.5056 & 9.4678E-03 & 1.5023
              & 9.9521E-03 & 1.5059 & 9.1113E-03 & 1.5003 & 9.5782E-03 & 1.5043 \\
\hline
1.99 & 64 & 1.3122E-01 & -- & 1.3757E-01 & -- & 1.3117E-01 & -- & 1.3752E-01 & -- & 1.3105E-01 & --
               & 1.3739E-01 & -- & 1.3119E-01 & -- & 1.3753E-01 & -- \\
     & 128 & 7.1292E-02 & 0.8802 & 7.4522E-02 & 0.8844 & 7.1318E-02 & 0.8791 & 7.4540E-02 & 0.8836 & 7.1297E-02 & 0.8782
               & 7.4508E-02 & 0.8828 & 7.1453E-02 & 0.8766 & 7.4661E-02 & 0.8813 \\
     & 256 & 3.5500E-02 & 1.0059 & 3.6809E-02 & 1.0176 & 3.5495E-02 & 1.0067 & 3.6804E-02 & 1.0182 & 3.5471E-02 & 1.0072
               & 3.6779E-02 & 1.0185 & 3.5526E-02 & 1.0081 & 3.6833E-02 & 1.0194 \\
     & 512 & 1.2807E-02 & 1.4709 & 1.3271E-02 & 1.4718 & 1.2809E-02 & 1.4705 & 1.3271E-02 & 1.4716 & 1.2803E-02 & 1.4702
               & 1.3265E-02 & 1.4713 & 1.2828E-02 & 1.4696 & 1.3289E-02 & 1.4708 \\
\hline
		\end{tabular}
		\label{tab1}
	\end{center}
\end{table}
\begin{table}[t]\scriptsize\tabcolsep=2.0pt
	\begin{center}
		\caption{The maximum norm errors and convergence orders for Example 1 with $\tau = h$.}
		\centering
		\begin{tabular}{cccccccccccccccccc}
			\hline
			& & \multicolumn{4}{c}{$\lambda = 0$} & \multicolumn{4}{c}{$\lambda = 1$}
			& \multicolumn{4}{c}{$\lambda = 5$} & \multicolumn{4}{c}{$\lambda = 10$} \\
			[-2pt] \cmidrule(lr){3-6} \cmidrule(lr){7-10} \cmidrule(lr){11-14} \cmidrule(lr){15-18}\\ [-11pt]
			& & \multicolumn{2}{c}{\rm{L-IES \eqref{eq2.4}}} & \multicolumn{2}{c}{\rm{NL-IES \eqref{eq2.3}}}
			& \multicolumn{2}{c}{\rm{L-IES \eqref{eq2.4}}} & \multicolumn{2}{c}{\rm{NL-IES \eqref{eq2.3}}}
			& \multicolumn{2}{c}{\rm{L-IES \eqref{eq2.4}}} & \multicolumn{2}{c}{\rm{NL-IES \eqref{eq2.3}}}
			& \multicolumn{2}{c}{\rm{L-IES \eqref{eq2.4}}} & \multicolumn{2}{c}{\rm{NL-IES \eqref{eq2.3}}} \\
			[-2pt] \cmidrule(lr){3-4} \cmidrule(lr){5-6} \cmidrule(lr){7-8} \cmidrule(lr){9-10}
			\cmidrule(lr){11-12} \cmidrule(lr){13-14} \cmidrule(lr){15-16} \cmidrule(lr){17-18} \\ [-11pt]
			$\alpha$ & $N$ & $Err(\tau,h)$ & $Order2$ & $Err(\tau,h)$ & $Order2$ & $Err(\tau,h)$ & $Order2$
			& $Err(\tau,h)$ & $Order2$ & $Err(\tau,h)$ & $Order2$ & $Err(\tau,h)$ & $Order2$
			& $Err(\tau,h)$ & $Order2$ & $Err(\tau,h)$ & $Order2$\\
\hline
1.1 & 64 & 2.4026E-01 & -- & 2.3738E-01 & -- & 2.8807E-01 & -- & 2.8824E-01 & -- & 4.8922E-01 & --
              & 4.9031E-01 & -- & 8.3296E-01 & -- & 8.4089E-01 & -- \\
    & 128 & 1.9494E-01 & 0.3016 & 1.9416E-01 & 0.2900 & 1.7183E-01 & 0.7454 & 1.7102E-01 & 0.7531 & 2.2356E-01 & 1.1298
              & 2.2247E-01 & 1.1401 & 3.8624E-01 & 1.1087 & 3.8625E-01 & 1.1224 \\
    & 256 & 1.5375E-01 & 0.3424 & 1.5350E-01 & 0.3390 & 1.0246E-01 & 0.7459 & 1.0222E-01 & 0.7425 & 9.4040E-02 & 1.2493
              & 9.3496E-02 & 1.2506 & 1.6193E-01 & 1.2541 & 1.6150E-01 & 1.2580 \\
    & 512 & 9.6958E-02 & 0.6652 & 9.6904E-02 & 0.6636 & 4.7287E-02 & 1.1155 & 4.7236E-02 & 1.1137 & 3.1677E-02 & 1.5698
              & 3.1520E-02 & 1.5686 & 5.4094E-02 & 1.5818 & 5.3945E-02 & 1.5820 \\
\hline
1.5 & 64 & 7.3489E-02 & -- & 7.2202E-02 & -- & 6.0097E-02 & -- & 5.9292E-02 & -- & 1.3157E-01 & --
              & 1.3375E-01 & -- & 3.0816E-01 & -- & 3.1344E-01 & -- \\
    & 128 & 4.8998E-02 & 0.5848 & 4.8626E-02 & 0.5703 & 3.7724E-02 & 0.6718 & 3.7582E-02 & 0.6578 & 6.7287E-02 & 0.9674
              & 6.8079E-02 & 0.9743 & 1.7000E-01 & 0.8581 & 1.7235E-01 & 0.8628 \\
    & 256 & 2.5042E-02 & 0.9684 & 2.4956E-02 & 0.9623 & 1.9558E-02 & 0.9477 & 1.9515E-02 & 0.9455 & 3.0301E-02 & 1.1510
              & 3.0576E-02 & 1.1548 & 8.0107E-02 & 1.0855 & 8.0929E-02 & 1.0906 \\
    & 512 & 9.2783E-03 & 1.4324 & 9.2628E-03 & 1.4299 & 7.5546E-03 & 1.3723 & 7.5455E-03 & 1.3709 & 1.0369E-02 & 1.5471
              & 1.0449E-02 & 1.5490 & 2.7766E-02 & 1.5286 & 2.7995E-02 & 1.5315 \\
\hline
1.9 & 64 & 1.1145E-01 & -- & 1.1829E-01 & -- & 1.0157E-01 & -- & 1.0836E-01 & -- & 9.6149E-02 & --
              & 1.0316E-01 & -- & 1.0809E-01 & -- & 1.1351E-01 & -- \\
    & 128 & 6.0209E-02 & 0.8883 & 6.3279E-02 & 0.9025 & 5.4910E-02 & 0.8873 & 5.7978E-02 & 0.9023 & 5.2520E-02 & 0.8724
              & 5.5832E-02 & 0.8857 & 6.0055E-02 & 0.8479 & 6.3014E-02 & 0.8491 \\
    & 256 & 2.8477E-02 & 1.0802 & 2.9948E-02 & 1.0793 & 2.6153E-02 & 1.0701 & 2.7608E-02 & 1.0704 & 2.4581E-02 & 1.0953
              & 2.6008E-02 & 1.1021 & 2.8130E-02 & 1.0942 & 2.9453E-02 & 1.0973 \\
    & 512 & 1.0045E-02 & 1.5033 & 1.0538E-02 & 1.5069 & 9.2257E-03 & 1.5032 & 9.7132E-03 & 1.5071 & 8.5957E-03 & 1.5159
              & 9.0714E-03 & 1.5196 & 9.8009E-03 & 1.5211 & 1.0259E-02 & 1.5215 \\
\hline
1.99 & 64 & 1.3155E-01 & -- & 1.3788E-01 & -- & 1.2545E-01 & -- & 1.3173E-01 & -- & 1.2630E-01 & --
               & 1.3326E-01 & -- & 1.3911E-01 & -- & 1.4597E-01 & -- \\
     & 128 & 7.1323E-02 & 0.8832 & 7.4580E-02 & 0.8866 & 6.8822E-02 & 0.8662 & 7.2074E-02 & 0.8700 & 6.8825E-02 & 0.8759
               & 7.1952E-02 & 0.8891 & 7.7047E-02 & 0.8524 & 8.0114E-02 & 0.8655 \\
     & 256 & 3.5498E-02 & 1.0066 & 3.6803E-02 & 1.0190 & 3.4089E-02 & 1.0136 & 3.5389E-02 & 1.0262 & 3.2523E-02 & 1.0815
               & 3.3934E-02 & 1.0843 & 3.6251E-02 & 1.0877 & 3.7635E-02 & 1.0900 \\
     & 512 & 1.2802E-02 & 1.4714 & 1.3265E-02 & 1.4722 & 1.2329E-02 & 1.4673 & 1.2790E-02 & 1.4683 & 1.1514E-02 & 1.4981
               & 1.1989E-02 & 1.5010 & 1.2871E-02 & 1.4939 & 1.3329E-02 & 1.4975 \\
\hline
		\end{tabular}
		\label{tab2}
	\end{center}
\end{table}
\end{landscape}
\begin{table}[t]\footnotesize\tabcolsep=2.0pt
	\begin{center}
		\caption{The maximum norm errors between NL-IES \eqref{eq2.3} and L-IES \eqref{eq2.4} for Example 1
			with $h = 2^{-10}$.}
		\centering
		\begin{tabular}{cccccccccc}
			\hline
			& & \multicolumn{2}{c}{$\lambda = 0$} & \multicolumn{2}{c}{$\lambda = 1$}
			& \multicolumn{2}{c}{$\lambda = 5$} & \multicolumn{2}{c}{$\lambda = 10$} \\
			[-2pt] \cmidrule(lr){3-4} \cmidrule(lr){5-6} \cmidrule(lr){7-8} \cmidrule(lr){9-10}\\ [-11pt]
			$\alpha$ & $M$ & $Err(\tau,h)$ & $Order1$ & $Err(\tau,h)$ & $Order1$ & $Err(\tau,h)$ & $Order1$
			& $Err(\tau,h)$ & $Order1$ \\
\hline
1.1 & 64 & 1.1812E-02 & -- & 7.9283E-03 & -- & 6.5488E-03 & -- & 5.8862E-03 & -- \\
    & 128 & 6.5390E-03 & 0.8531 & 3.9934E-03 & 0.9894 & 3.2863E-03 & 0.9948 & 2.9518E-03 & 0.9957 \\
    & 256 & 3.4857E-03 & 0.9076 & 2.0043E-03 & 0.9945 & 1.6467E-03 & 0.9969 & 1.4784E-03 & 0.9976 \\
    & 512 & 1.8064E-03 & 0.9483 & 1.0041E-03 & 0.9972 & 8.2423E-04 & 0.9985 & 7.3988E-04 & 0.9987 \\
\hline
1.5 & 64 & 8.9765E-03 & -- & 8.8547E-03 & -- & 8.6193E-03 & -- & 8.5281E-03 & -- \\
    & 128 & 4.6220E-03 & 0.9576 & 4.5559E-03 & 0.9587 & 4.4082E-03 & 0.9674 & 4.3454E-03 & 0.9727 \\
    & 256 & 2.4014E-03 & 0.9446 & 2.3235E-03 & 0.9714 & 2.2390E-03 & 0.9773 & 2.1958E-03 & 0.9847 \\
    & 512 & 1.2244E-03 & 0.9718 & 1.1764E-03 & 0.9819 & 1.1280E-03 & 0.9891 & 1.1049E-03 & 0.9908 \\
\hline
1.9 & 64 & 8.1415E-03 & -- & 8.1710E-03 & -- & 8.2299E-03 & -- & 8.2611E-03 & -- \\
    & 128 & 4.1759E-03 & 0.9632 & 4.1567E-03 & 0.9751 & 4.2044E-03 & 0.9690 & 4.2457E-03 & 0.9603 \\
    & 256 & 2.2051E-03 & 0.9212 & 2.2097E-03 & 0.9116 & 2.2166E-03 & 0.9236 & 2.2180E-03 & 0.9367 \\
    & 512 & 1.1399E-03 & 0.9519 & 1.1377E-03 & 0.9577 & 1.1312E-03 & 0.9705 & 1.1329E-03 & 0.9692 \\
\hline
1.99 & 64 & 7.7831E-03 & -- & 7.7909E-03 & -- & 7.8126E-03 & -- & 7.8311E-03 & -- \\
     & 128 & 4.1471E-03 & 0.9082 & 4.1464E-03 & 0.9099 & 4.1452E-03 & 0.9144 & 4.1443E-03 & 0.9181 \\
     & 256 & 2.1825E-03 & 0.9261 & 2.1816E-03 & 0.9265 & 2.1800E-03 & 0.9271 & 2.1790E-03 & 0.9275 \\
     & 512 & 1.1267E-03 & 0.9539 & 1.1264E-03 & 0.9537 & 1.1260E-03 & 0.9531 & 1.1258E-03 & 0.9527 \\	
\hline
\end{tabular}
\label{tab3}
\end{center}
\end{table}	
\begin{table}[H]\footnotesize\tabcolsep=2.0pt
	\begin{center}
		\caption{The maximum norm errors between NL-IES \eqref{eq2.3} and L-IES \eqref{eq2.4} for Example 1
			with $\tau = h$.}
		\centering
		\begin{tabular}{cccccccccc}
			\hline
			& & \multicolumn{2}{c}{$\lambda = 0$} & \multicolumn{2}{c}{$\lambda = 1$}
			& \multicolumn{2}{c}{$\lambda = 5$} & \multicolumn{2}{c}{$\lambda = 10$} \\
			[-2pt] \cmidrule(lr){3-4} \cmidrule(lr){5-6} \cmidrule(lr){7-8} \cmidrule(lr){9-10} \\ [-11pt]
			$\alpha$ & $N$ & $Err(\tau,h)$ & $Order1$ & $Err(\tau,h)$ & $Order1$ & $Err(\tau,h)$ & $Order1$
			& $Err(\tau,h)$ & $Order1$ \\
			\hline
			1.1 & 64 & 1.0573E-02 & -- & 8.3127E-03 & -- & 7.9461E-03 & -- & 8.1361E-03 & -- \\
			& 128 & 6.0969E-03 & 0.7942 & 4.0972E-03 & 1.0207 & 3.7134E-03 & 1.0975 & 3.6526E-03 & 1.1554 \\
			& 256 & 3.3629E-03 & 0.8584 & 2.0262E-03 & 1.0159 & 1.7471E-03 & 1.0878 & 1.6332E-03 & 1.1612 \\
			& 512 & 1.7831E-03 & 0.9153 & 1.0076E-03 & 1.0079 & 8.4131E-04 & 1.0543 & 7.6363E-04 & 1.0968 \\
			\hline
			1.5 & 64 & 8.9043E-03 & -- & 8.8174E-03 & -- & 8.5584E-03 & -- & 8.4458E-03 & -- \\
			& 128 & 4.6102E-03 & 0.9497 & 4.5477E-03 & 0.9552 & 4.4088E-03 & 0.9570 & 4.3427E-03 & 0.9596 \\
			& 256 & 2.3971E-03 & 0.9435 & 2.3200E-03 & 0.9710 & 2.2381E-03 & 0.9781 & 2.2016E-03 & 0.9800 \\
			& 512 & 1.2235E-03 & 0.9703 & 1.1758E-03 & 0.9805 & 1.1287E-03 & 0.9876 & 1.1061E-03 & 0.9931 \\
			\hline
			1.9 & 64 & 8.1434E-03 & -- & 8.1893E-03 & -- & 8.2678E-03 & -- & 8.2806E-03 & -- \\
			& 128 & 4.1773E-03 & 0.9631 & 4.1606E-03 & 0.9769 & 4.2237E-03 & 0.9690 & 4.2707E-03 & 0.9553 \\
			& 256 & 2.2049E-03 & 0.9219 & 2.2099E-03 & 0.9128 & 2.2176E-03 & 0.9295 & 2.2178E-03 & 0.9453 \\
			& 512 & 1.1399E-03 & 0.9518 & 1.1378E-03 & 0.9577 & 1.1310E-03 & 0.9714 & 1.1330E-03 & 0.9690 \\
			\hline
			1.99 & 64 & 7.7944E-03 & -- & 7.8275E-03 & -- & 7.9290E-03 & -- & 8.0097E-03 & -- \\
			& 128 & 4.1447E-03 & 0.9112 & 4.1426E-03 & 0.9180 & 4.1333E-03 & 0.9398 & 4.1162E-03 & 0.9604 \\
			& 256 & 2.1827E-03 & 0.9252 & 2.1817E-03 & 0.9251 & 2.1795E-03 & 0.9233 & 2.1772E-03 & 0.9188 \\
			& 512 & 1.1267E-03 & 0.9540 & 1.1264E-03 & 0.9537 & 1.1260E-03 & 0.9528 & 1.1257E-03 & 0.9517 \\
			\hline
		\end{tabular}
		\label{tab4}
	\end{center}
\end{table}
Tables \ref{tab1}-\ref{tab2} show that the convergence orders
of the implicit schemes NL-IES \eqref{eq2.3} and L-IES \eqref{eq2.4}
for different $\alpha$ and $\lambda$ can indeed reach $1$ in both time and space.
The errors between NL-IES \eqref{eq2.3} and L-IES \eqref{eq2.4} are listed in Tables \ref{tab3}-\ref{tab4}.
From such tables, the $Order1$ and $Order2$ are almost 1,
which in accord with the result at the end of Section \ref{sec3}.
In Table \ref{tab5}, the CPU time and number of iterations of the methods BS,
$\mathcal{I}$ and $\mathcal{P}$ are reported.
The method $\mathcal{I}$ in most cases needs more than 1000 iterative steps
to obtain the solutions of $J^k \bm{z} = -\bm{f}(\bm{u}^{(k)})$,
which implies that the Jacobian matrices are very ill-conditioned.
For the method $\mathcal{P}$, the $\mathrm{Iter2}$ is greatly reduced compared with the method $\mathcal{I}$.
This means that our preconditioner $P_{\ell}$ is efficient for solving the Jacobian equations in Algorithm \ref{alg1},
but the $\mathrm{Iter2}$ grows slightly fast in several cases such as $(\alpha, \lambda) = (1.9,5)$.
On the other hand, as seen from Table \ref{tab5}, the total CPU time
of the method $\mathcal{P}$ is the smallest one among them.
The eigenvalues of the initial Jacobian matrix $J^{0}$
and its preconditioned matrix $P_{\ell}^{-1} J^{0}$ are drawn in Fig. \ref{fig2}.
As can be seen, the eigenvalues of $P_{\ell}^{-1} J^{0}$ are clustered around $1$.

\noindent \textbf{Example 2.} In this example, we consider Eq. \eqref{eq1.1} with $T = 1$, the initial value
$u(x,0) = \frac{4 \exp(10 x)}{\left(\exp(10 x) + 1 \right)^2} ~(x \in [-1,1])$,
the nonlinear source term $f(u(x,t),x,t) = -u(x,t) \left(1 - u(x,t)\right)$
and the discontinuous coefficients
\begin{equation*}
d_{+}(x) =
\begin{cases}
1.5 \exp(-x), ~ -1 \leq x < 0, \\
2 \mathrm{sech}(x), ~ 0 \leq x \leq 1,
\end{cases}\quad
d_{-}(x) =
\begin{cases}
\exp(x), ~ -1 \leq x < 0, \\
0.1 + \mathrm{sech}(-x), ~ 0 \leq x \leq 1.
\end{cases}
\end{equation*}
Similar to Example 1, we regard the numerical solution on the finer mesh ($M = N = 1024$) as our exact solution.

\begin{figure}[t]
	\centering
	\subfigure[Eigenvalues of $J^{0}$]
	{\includegraphics[width=3.0in,height=1.4in]{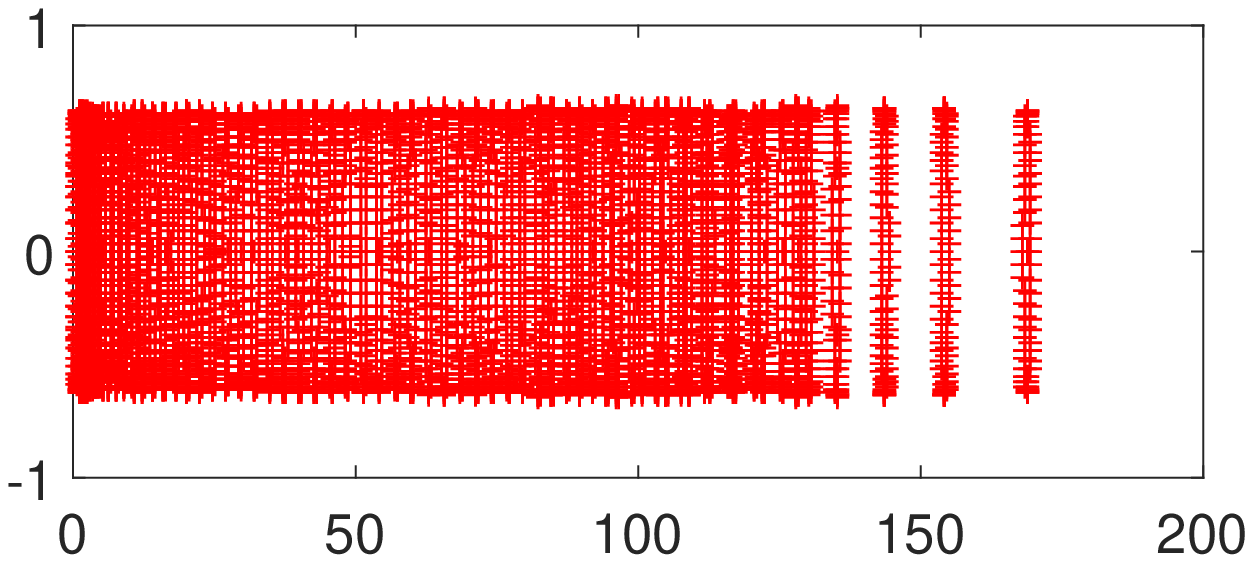}}\hspace{0.2in}
	\subfigure[Eigenvalues of $P_{\ell}^{-1} J^{0}$]
	{\includegraphics[width=3.0in,height=1.4in]{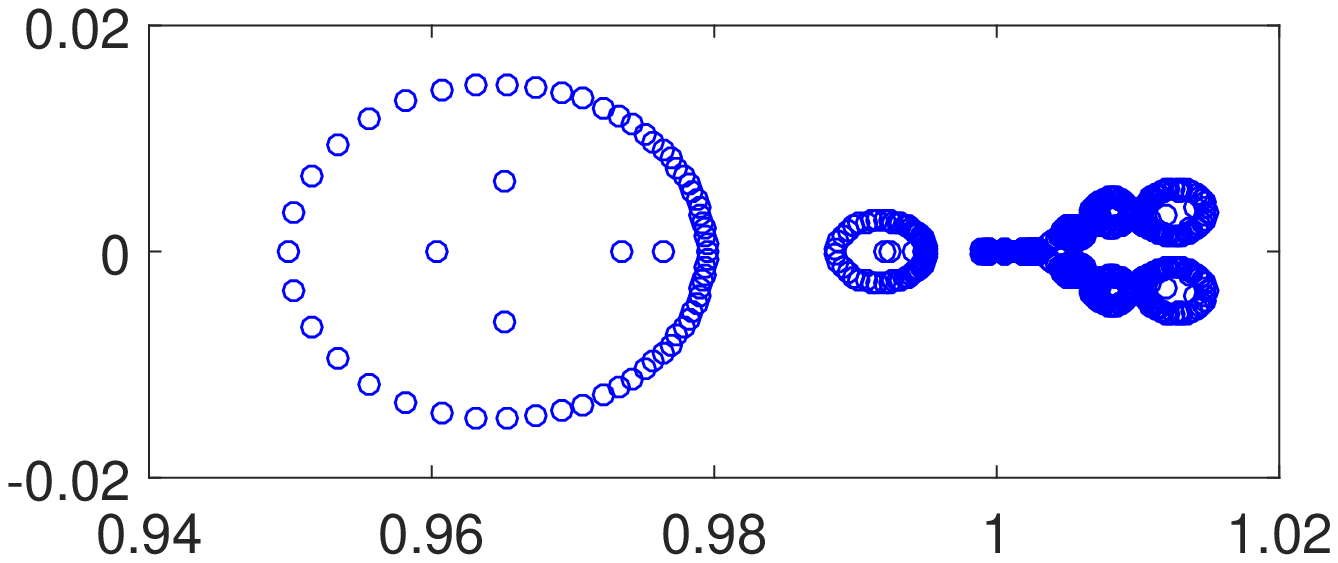}}\vspace{4em}
	\subfigure[Eigenvalues of $J^{0}$]
	{\includegraphics[width=3.0in,height=1.4in]{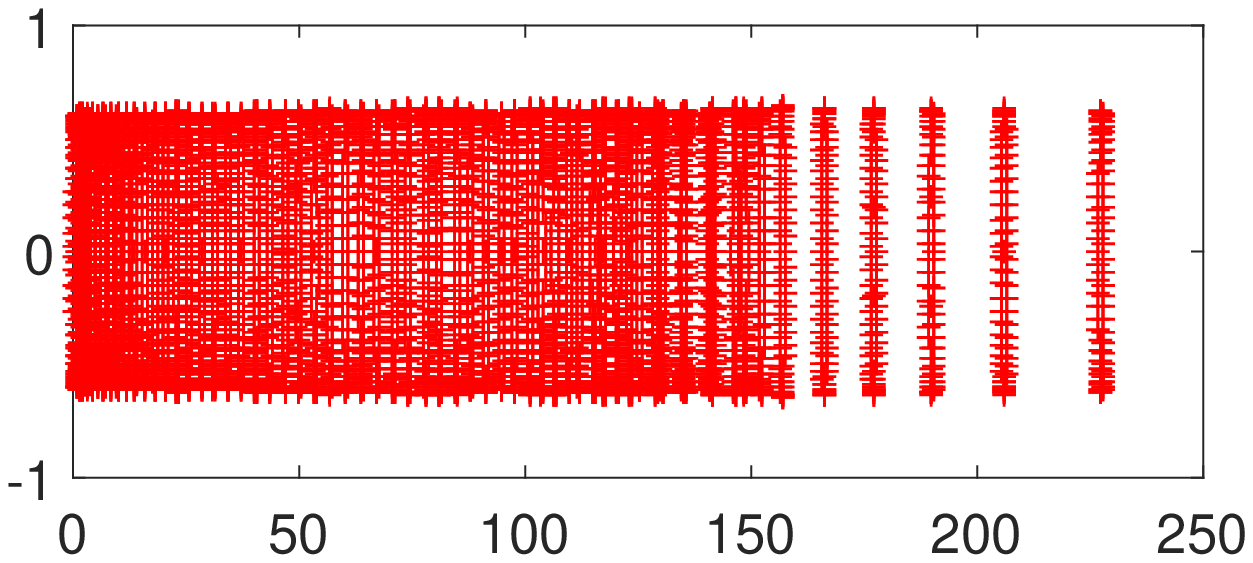}}\hspace{0.2in}
	\subfigure[Eigenvalues of $P_{\ell}^{-1} J^{0}$]
	{\includegraphics[width=3.0in,height=1.4in]{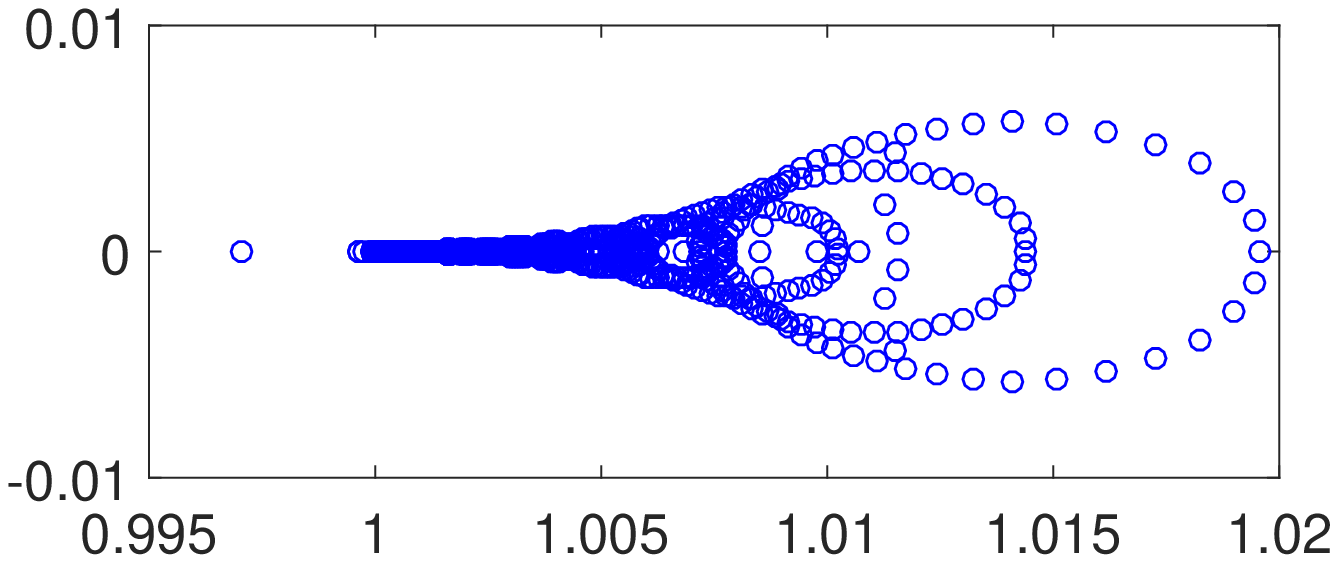}}
	\caption{Spectra of $J^{0}$ and $P_{\ell}^{-1} J^{0}$, when $\alpha = 1.9$, $M = N = 65$ in Example 2.
		Top row: $\lambda = 0$; Bottom row: $\lambda = 10$.}
	\label{fig4}
\end{figure}

It can be seen from Tables \ref{tab6}-\ref{tab7} that the convergence orders of the two implicit schemes
are $1$ in both time and space for the discontinuous coefficients.
Tables \ref{tab8} and \ref{tab9} display the errors between NL-IES \eqref{eq2.3} and L-IES \eqref{eq2.4},
and the rates of such errors can indeed reach $1$.
The performance of the method $\mathcal{P}$ shown in Table \ref{tab10} is the best one among them
in aspects of CPU time and the number of iterations.
The $\mathrm{Iter1}$ of the methods BS, $\mathcal{I}$ and $\mathcal{P}$ is small,
which indicates that the initial vector $\bm{u}^{(0)}$ is a good enough initial value.
As for the $\mathrm{Iter2}$, the method $\mathcal{P}$ requires less iterative steps
than the method $\mathcal{I}$ under the same termination condition.
This illustrates that our proposed preconditioner $P_{\ell}$ is efficient
and can accelerate solving the Jacobian equations in Algorithm \ref{alg1}.
Furthermore, Fig. \ref{fig4} displays the eigenvalues of $J^{0}$ and $P_{\ell}^{-1} J^{0}$.
\begin{table}[H]\footnotesize\tabcolsep=3.0pt
	\begin{center}
		\caption{Results of different methods when $M = N$ for Example 1.}
		\centering
		\begin{tabular}{cccccccc}
			\hline
			& & \multicolumn{2}{c}{\rm{BS}} & \multicolumn{2}{c}{$\mathcal{I}$} & \multicolumn{2}{c}{$\mathcal{P}$} \\
			[-2pt] \cmidrule(lr){3-4} \cmidrule(lr){5-6} \cmidrule(lr){7-8} \\ [-11pt]
			$(\alpha, \lambda)$ & $N$ & $\mathrm{Iter1}$ & $\textrm{Time}$ & $(\mathrm{Iter1}, \mathrm{Iter2})$ & $\textrm{Time}$
			& $(\mathrm{Iter1}, \mathrm{Iter2})$ & $\textrm{Time}$ \\
			\hline
			(1.1, 0) & 129 & 5.0 & 0.509 & (5.0, 230.4) & 5.980 & (5.0, 4.8) & 0.199 \\
			         & 257 & 5.0 & 3.997 & (5.0, 521.0) & 53.085 & (5.0, 6.6) & 0.928 \\
			         & 513 & 5.0 & 35.110 & (6.0, 688.8) & 315.968 & (5.0, 9.8) & 4.704 \\
			         & 1025 & \dag & \dag & \ddag & \ddag & (5.0, 16.2) & 26.179 \\
			\hline
			(1.5, 0) & 129 & 4.0 & 0.448 & (4.0, 483.3) & 10.255 & (4.0, 12.8) & 0.399 \\
			         & 257 & 4.0 & 3.336 & (5.0, 817.4) & 84.057 & (4.0, 25.5) & 2.380 \\
			         & 513 & 4.0 & 29.445 & \ddag & \ddag & (4.0, 62.3) & 22.153 \\
			         & 1025 & \dag & \dag & \ddag & \ddag & (4.0, 169.5) & 207.366 \\
			\hline
			(1.9, 0) & 129 & 4.0 & 0.441 & (5.0, 913.2) & 27.381 & (4.0, 9.3) & 0.290 \\
			         & 257 & 4.0 & 3.393 & \ddag & \ddag & (4.0, 22.3) & 2.248 \\
			         & 513 & 4.0 & 29.532 & \ddag & \ddag & (4.0, 66.5) & 23.564 \\
			         & 1025 & \dag & \dag & \ddag & \ddag & (4.0, 215.0) & 262.269 \\
			\hline
			(1.1, 5) & 129 & 6.0 & 0.595 & (6.0, 252.7) & 9.361 & (6.0, 3.8) & 0.178 \\
			         & 257 & 6.0 & 4.666 & (24.0, 332.5) & 161.449 & (6.0, 4.5) & 0.760 \\
			         & 513 & 6.0 & 41.446 & \ddag & \ddag & (6.0, 7.0) & 4.022 \\
			         & 1025 & \dag & \dag & \ddag & \ddag & (7.0, 13.4) & 29.704 \\
			\hline
			(1.5, 5) & 129 & 4.0 & 0.454 & (5.0, 628.8) & 19.502 & (5.0, 7.6) & 0.283 \\
			         & 257 & 5.0 & 4.002 & \ddag & \ddag & (5.0, 18.8) & 2.381 \\
			         & 513 & 5.0 & 35.842 & \ddag & \ddag & (5.0, 52.2) & 23.144 \\
			         & 1025 & \dag & \dag & \ddag & \ddag & (5.0, 150.4) & 227.671 \\
			\hline
			(1.9, 5) & 129 & 4.0 & 0.457 & (5.0, 944.4) & 28.843 & (4.0, 6.3) & 0.190 \\
			         & 257 & 4.0 & 3.350 & \ddag & \ddag & (4.0, 17.0) & 1.685 \\
			         & 513 & 4.0 & 29.932 & \ddag & \ddag & (4.0, 56.3) & 19.871 \\
			         & 1025 & \dag & \dag & \ddag & \ddag & (4.0, 218.3) & 263.273 \\
			\hline
			(1.1, 10) & 129 & 6.0 & 0.691 & (6.0, 242.8) & 9.137 & (6.0, 3.7) & 0.171 \\
			          & 257 & 7.0 & 5.341 & (7.0, 519.9) & 74.813 & (7.0, 4.3) & 0.727 \\
			          & 513 & 7.0 & 47.171 & \ddag & \ddag & (7.0, 6.1) & 4.011 \\
			          & 1025 & \dag & \dag & \ddag & \ddag & (7.0, 11.7) & 26.043 \\
			\hline
			(1.5, 10) & 129 & 5.0 & 0.542 & (5.0, 664.4) & 20.579 & (5.0, 5.4) & 0.221 \\
			          & 257 & 5.0 & 4.007 & (19.0, 373.5) & 146.726 & (5.0, 14.4) & 1.820 \\
			          & 513 & 5.0 & 36.374 & \ddag & \ddag & (5.0, 43.6) & 19.333 \\
			          & 1025 & \dag & \dag & \ddag & \ddag & (5.0, 143.2) & 215.481 \\
			\hline
			(1.9, 10) & 129 & 4.0 & 0.484 & (5.0, 898.4) & 28.227 & (4.0, 4.8) & 0.153 \\
			          & 257 & 4.0 & 3.327 & \ddag & \ddag & (4.0, 12.3) & 1.253 \\
			          & 513 & 4.0 & 30.498 & \ddag & \ddag & (4.0, 51.3) & 18.470 \\
			          & 1025 & \dag & \dag & \ddag & \ddag & (4.0, 186.0) & 225.510 \\ 			
			\hline
\end{tabular}
\label{tab5}
\end{center}
\end{table}
\begin{landscape}
	\begin{table}[t]\scriptsize\tabcolsep=2.0pt
		\begin{center}
			\caption{The maximum norm errors and convergence orders for Example 2 with $h = 2^{-10}$.}
			\centering
			\begin{tabular}{cccccccccccccccccc}
				\hline
				& & \multicolumn{4}{c}{$\lambda = 0$} & \multicolumn{4}{c}{$\lambda = 1$}
				& \multicolumn{4}{c}{$\lambda = 5$} & \multicolumn{4}{c}{$\lambda = 10$} \\
				[-2pt] \cmidrule(lr){3-6} \cmidrule(lr){7-10} \cmidrule(lr){11-14} \cmidrule(lr){15-18}\\ [-11pt]
				& & \multicolumn{2}{c}{\rm{L-IES \eqref{eq2.4}}} & \multicolumn{2}{c}{\rm{NL-IES \eqref{eq2.3}}}
				& \multicolumn{2}{c}{\rm{L-IES \eqref{eq2.4}}} & \multicolumn{2}{c}{\rm{NL-IES \eqref{eq2.3}}}
				& \multicolumn{2}{c}{\rm{L-IES \eqref{eq2.4}}} & \multicolumn{2}{c}{\rm{NL-IES \eqref{eq2.3}}}
				& \multicolumn{2}{c}{\rm{L-IES \eqref{eq2.4}}} & \multicolumn{2}{c}{\rm{NL-IES \eqref{eq2.3}}} \\
				[-2pt] \cmidrule(lr){3-4} \cmidrule(lr){5-6} \cmidrule(lr){7-8} \cmidrule(lr){9-10}
				\cmidrule(lr){11-12} \cmidrule(lr){13-14} \cmidrule(lr){15-16} \cmidrule(lr){17-18} \\ [-11pt]
				$\alpha$ & $M$ & $Err(\tau,h)$ & $Order1$ & $Err(\tau,h)$ & $Order1$ & $Err(\tau,h)$ & $Order1$
				& $Err(\tau,h)$ & $Order1$ & $Err(\tau,h)$ & $Order1$ & $Err(\tau,h)$ & $Order1$
				& $Err(\tau,h)$ & $Order1$ & $Err(\tau,h)$ & $Order1$\\
				\hline
				1.1 & 64 & 3.2431E-02 & -- & 3.3730E-02 & -- & 6.8046E-03 & -- & 4.9594E-03 & -- & 4.3160E-03 & --
				& 2.3456E-03 & -- & 3.4101E-03 & -- & 1.3524E-03 & -- \\
				& 128 & 1.6340E-02 & 0.9890 & 1.6918E-02 & 0.9955 & 3.2092E-03 & 1.0843 & 2.3436E-03 & 1.0814 & 2.0253E-03 & 1.0916
				& 1.1032E-03 & 1.0883 & 1.5961E-03 & 1.0953 & 6.3392E-04 & 1.0931 \\
				& 256 & 7.3028E-03 & 1.1619 & 7.5422E-03 & 1.1655 & 1.3830E-03 & 1.2144 & 1.0110E-03 & 1.2129 & 8.7035E-04 & 1.2185
				& 4.7469E-04 & 1.2166 & 6.8508E-04 & 1.2202 & 2.7231E-04 & 1.2191 \\
				& 512 & 2.4886E-03 & 1.5531 & 2.5666E-03 & 1.5551 & 4.6229E-04 & 1.5809 & 3.3812E-04 & 1.5802 & 2.9052E-04 & 1.5830
				& 1.5855E-04 & 1.5820 & 2.2853E-04 & 1.5839 & 9.0900E-05 & 1.5829 \\
				\hline
				1.5 & 64 & 5.2358E-02 & -- & 5.0868E-02 & -- & 4.3947E-02 & -- & 4.2489E-02 & -- & 3.1695E-02 & --
				& 3.0122E-02 & -- & 2.6224E-02 & -- & 2.4699E-02 & -- \\
				& 128 & 2.7642E-02 & 0.9215 & 2.6909E-02 & 0.9187 & 2.3157E-02 & 0.9243 & 2.2433E-02 & 0.9215 & 1.6578E-02 & 0.9350
				& 1.5849E-02 & 0.9264 & 1.3353E-02 & 0.9737 & 1.2661E-02 & 0.9641 \\
				& 256 & 1.2690E-02 & 1.1232 & 1.2368E-02 & 1.1215 & 1.0624E-02 & 1.1241 & 1.0302E-02 & 1.1227 & 7.5608E-03 & 1.1327
				& 7.2512E-03 & 1.1281 & 5.9884E-03 & 1.1569 & 5.6980E-03 & 1.1519 \\
				& 512 & 4.4268E-03 & 1.5194 & 4.3191E-03 & 1.5178 & 3.6814E-03 & 1.5290 & 3.5748E-03 & 1.5270 & 2.6047E-03 & 1.5374
				& 2.5021E-03 & 1.5351 & 2.0525E-03 & 1.5448 & 1.9541E-03 & 1.5440 \\
				\hline
				1.9 & 64 & 8.5313E-02 & -- & 8.3849E-02 & -- & 8.3615E-02 & -- & 8.2142E-02 & -- & 8.0497E-02 & --
				& 7.9014E-02 & -- & 7.8786E-02 & -- & 7.7299E-02 & -- \\
				& 128 & 5.8082E-02 & 0.5547 & 5.7382E-02 & 0.5472 & 5.6620E-02 & 0.5625 & 5.5923E-02 & 0.5547 & 5.3634E-02 & 0.5858
				& 5.2945E-02 & 0.5776 & 5.1839E-02 & 0.6039 & 5.1155E-02 & 0.5956 \\
				& 256 & 2.9567E-02 & 0.9741 & 2.9299E-02 & 0.9697 & 2.8709E-02 & 0.9798 & 2.8444E-02 & 0.9753 & 2.6837E-02 & 0.9989
				& 2.6580E-02 & 0.9942 & 2.5645E-02 & 1.0154 & 2.5392E-02 & 1.0105 \\
				& 512 & 1.1116E-02 & 1.4114 & 1.1021E-02 & 1.4106 & 1.0791E-02 & 1.4117 & 1.0697E-02 & 1.4109 & 1.0083E-02 & 1.4123
				& 9.9912E-03 & 1.4116 & 9.6308E-03 & 1.4130 & 9.5406E-03 & 1.4122 \\
				\hline
				1.99 & 64 & 8.7279E-02 & -- & 8.5852E-02 & -- & 8.7128E-02 & -- & 8.5699E-02 & -- & 8.6886E-02 & --
				& 8.5455E-02 & -- & 8.6792E-02 & -- & 8.5360E-02 & -- \\
				& 128 & 6.2942E-02 & 0.4716 & 6.2243E-02 & 0.4639 & 6.2805E-02 & 0.4723 & 6.2106E-02 & 0.4645 & 6.2554E-02 & 0.4740
				& 6.1855E-02 & 0.4663 & 6.2449E-02 & 0.4749 & 6.1750E-02 & 0.4671 \\
				& 256 & 3.4261E-02 & 0.8775 & 3.3986E-02 & 0.8730 & 3.4176E-02 & 0.8779 & 3.3901E-02 & 0.8734 & 3.4008E-02 & 0.8792
				& 3.3734E-02 & 0.8747 & 3.3932E-02 & 0.8800 & 3.3658E-02 & 0.8755 \\
				& 512 & 1.2865E-02 & 1.4131 & 1.2768E-02 & 1.4124 & 1.2833E-02 & 1.4131 & 1.2736E-02 & 1.4124 & 1.2769E-02 & 1.4132
				& 1.2673E-02 & 1.4124 & 1.2741E-02 & 1.4132 & 1.2644E-02 & 1.4125 \\
				\hline
			\end{tabular}
			\label{tab6}
		\end{center}
	\end{table}
	\begin{table}[t]\scriptsize\tabcolsep=2.0pt
		\begin{center}
			\caption{The maximum norm errors and convergence orders for Example 2 with $\tau = h$.}
			\centering
			\begin{tabular}{cccccccccccccccccc}
				\hline
				& & \multicolumn{4}{c}{$\lambda = 0$} & \multicolumn{4}{c}{$\lambda = 1$}
				& \multicolumn{4}{c}{$\lambda = 5$} & \multicolumn{4}{c}{$\lambda = 10$} \\
				[-2pt] \cmidrule(lr){3-6} \cmidrule(lr){7-10} \cmidrule(lr){11-14} \cmidrule(lr){15-18}\\ [-11pt]
				& & \multicolumn{2}{c}{\rm{L-IES \eqref{eq2.4}}} & \multicolumn{2}{c}{\rm{NL-IES \eqref{eq2.3}}}
				& \multicolumn{2}{c}{\rm{L-IES \eqref{eq2.4}}} & \multicolumn{2}{c}{\rm{NL-IES \eqref{eq2.3}}}
				& \multicolumn{2}{c}{\rm{L-IES \eqref{eq2.4}}} & \multicolumn{2}{c}{\rm{NL-IES \eqref{eq2.3}}}
				& \multicolumn{2}{c}{\rm{L-IES \eqref{eq2.4}}} & \multicolumn{2}{c}{\rm{NL-IES \eqref{eq2.3}}} \\
				[-2pt] \cmidrule(lr){3-4} \cmidrule(lr){5-6} \cmidrule(lr){7-8} \cmidrule(lr){9-10}
				\cmidrule(lr){11-12} \cmidrule(lr){13-14} \cmidrule(lr){15-16} \cmidrule(lr){17-18} \\ [-11pt]
				$\alpha$ & $N$ & $Err(\tau,h)$ & $Order2$ & $Err(\tau,h)$ & $Order2$ & $Err(\tau,h)$ & $Order2$
				& $Err(\tau,h)$ & $Order2$ & $Err(\tau,h)$ & $Order2$ & $Err(\tau,h)$ & $Order2$
				& $Err(\tau,h)$ & $Order2$ & $Err(\tau,h)$ & $Order2$\\
				\hline
				1.1 & 64 & 1.0317E-01 & -- & 1.0452E-01 & -- & 1.3346E-01 & -- & 1.3487E-01 & -- & 2.1802E-01 & --
				& 2.1913E-01 & -- & 3.1633E-01 & -- & 3.1681E-01 & -- \\
				& 128 & 6.8449E-02 & 0.5919 & 6.8706E-02 & 0.6053 & 7.7164E-02 & 0.7904 & 7.7912E-02 & 0.7917 & 1.3089E-01 & 0.7361
				& 1.3152E-01 & 0.7365 & 1.9729E-01 & 0.6811 & 1.9769E-01 & 0.6804 \\
				& 256 & 4.5280E-02 & 0.5962 & 4.5395E-02 & 0.5979 & 3.8083E-02 & 1.0188 & 3.8428E-02 & 1.0197 & 6.6996E-02 & 0.9662
				& 6.7307E-02 & 0.9665 & 1.0488E-01 & 0.9116 & 1.0511E-01 & 0.9113 \\
				& 512 & 2.6552E-02 & 0.7701 & 2.6588E-02 & 0.7718 & 1.3841E-02 & 1.4602 & 1.3963E-02 & 1.4605 & 2.4969E-02 & 1.4239
				& 2.5081E-02 & 1.4242 & 4.0448E-02 & 1.3746 & 4.0540E-02 & 1.3745 \\
				\hline
				1.5 & 64 & 4.5283E-02 & -- & 4.3847E-02 & -- & 2.8941E-02 & -- & 2.7535E-02 & -- & 3.5242E-02 & --
				& 3.4077E-02 & -- & 9.9995E-02 & -- & 9.8520E-02 & -- \\
				& 128 & 2.3622E-02 & 0.9388 & 2.2923E-02 & 0.9357 & 1.4966E-02 & 0.9514 & 1.4277E-02 & 0.9476 & 1.8483E-02 & 0.9311
				& 1.7949E-02 & 0.9249 & 5.3729E-02 & 0.8962 & 5.3096E-02 & 0.8918 \\
				& 256 & 1.0920E-02 & 1.1132 & 1.0622E-02 & 1.1097 & 7.0537E-03 & 1.0852 & 6.7675E-03 & 1.0770 & 8.3827E-03 & 1.1407
				& 8.1585E-03 & 1.1375 & 2.4510E-02 & 1.1323 & 2.4265E-02 & 1.1297 \\
				& 512 & 3.8072E-03 & 1.5202 & 3.7073E-03 & 1.5186 & 2.4552E-03 & 1.5225 & 2.3625E-03 & 1.5183 & 2.8736E-03 & 1.5446
				& 2.8000E-03 & 1.5429 & 8.4061E-03 & 1.5439 & 8.3295E-03 & 1.5426 \\
				\hline
				1.9 & 64 & 8.2087E-02 & -- & 8.0640E-02 & -- & 7.8230E-02 & -- & 7.6775E-02 & -- & 6.5789E-02 & --
				& 6.4329E-02 & -- & 5.0179E-02 & -- & 4.8724E-02 & -- \\
				& 128 & 5.6383E-02 & 0.5419 & 5.5686E-02 & 0.5342 & 5.3996E-02 & 0.5349 & 5.3301E-02 & 0.5265 & 4.7243E-02 & 0.4777
				& 4.6554E-02 & 0.4666 & 4.0154E-02 & 0.3215 & 3.9467E-02 & 0.3040 \\
				& 256 & 2.8883E-02 & 0.9650 & 2.8616E-02 & 0.9605 & 2.7681E-02 & 0.9640 & 2.7417E-02 & 0.9591 & 2.4505E-02 & 0.9470
				& 2.4248E-02 & 0.9410 & 2.1569E-02 & 0.8966 & 2.1315E-02 & 0.8888 \\
				& 512 & 1.0852E-02 & 1.4123 & 1.0757E-02 & 1.4115 & 1.0404E-02 & 1.4118 & 1.0310E-02 & 1.4110 & 9.2250E-03 & 1.4095
				& 9.1333E-03 & 1.4087 & 8.1620E-03 & 1.4020 & 8.0715E-03 & 1.4010 \\
				\hline
				1.99 & 64 & 8.4400E-02 & -- & 8.2988E-02 & -- & 8.2930E-02 & -- & 8.1518E-02 & -- & 7.6184E-02 & --
				& 7.4777E-02 & -- & 6.5522E-02 & -- & 6.4128E-02 & -- \\
				& 128 & 6.1409E-02 & 0.4588 & 6.0712E-02 & 0.4509 & 6.0689E-02 & 0.4505 & 5.9992E-02 & 0.4423 & 5.7807E-02 & 0.3982
				& 5.7111E-02 & 0.3888 & 5.3735E-02 & 0.2861 & 5.3038E-02 & 0.2739 \\
				& 256 & 3.3581E-02 & 0.8708 & 3.3306E-02 & 0.8662 & 3.3280E-02 & 0.8668 & 3.3005E-02 & 0.8621 & 3.2186E-02 & 0.8448
				& 3.1911E-02 & 0.8397 & 3.0799E-02 & 0.8030 & 3.0524E-02 & 0.7971 \\
				& 512 & 1.2614E-02 & 1.4126 & 1.2517E-02 & 1.4119 & 1.2502E-02 & 1.4125 & 1.2405E-02 & 1.4118 & 1.2106E-02 & 1.4107
				& 1.2009E-02 & 1.4099 & 1.1626E-02 & 1.4055 & 1.1529E-02 & 1.4047 \\
				\hline
			\end{tabular}
			\label{tab7}
		\end{center}
	\end{table}
\end{landscape}
\begin{table}[t]\footnotesize\tabcolsep=2.0pt
	\begin{center}
		\caption{The maximum norm errors between NL-IES \eqref{eq2.3} and L-IES \eqref{eq2.4} for Example 2
			with $h = 2^{-10}$.}
		\centering
		\begin{tabular}{cccccccccc}
			\hline
			& & \multicolumn{2}{c}{$\lambda = 0$} & \multicolumn{2}{c}{$\lambda = 1$}
			& \multicolumn{2}{c}{$\lambda = 5$} & \multicolumn{2}{c}{$\lambda = 10$} \\
			[-2pt] \cmidrule(lr){3-4} \cmidrule(lr){5-6} \cmidrule(lr){7-8} \cmidrule(lr){9-10}\\ [-11pt]
			$\alpha$ & $M$ & $Err(\tau,h)$ & $Order1$ & $Err(\tau,h)$ & $Order1$ & $Err(\tau,h)$ & $Order1$
			& $Err(\tau,h)$ & $Order1$ \\
			\hline
1.1 & 64 & 2.8218E-03 & -- & 2.1570E-03 & -- & 2.1603E-03 & -- & 2.2519E-03 & -- \\
    & 128 & 1.5181E-03 & 0.8943 & 1.0877E-03 & 0.9877 & 1.0841E-03 & 0.9947 & 1.1289E-03 & 0.9962 \\
    & 256 & 7.9372E-04 & 0.9356 & 5.4617E-04 & 0.9939 & 5.4301E-04 & 0.9974 & 5.6515E-04 & 0.9982 \\
    & 512 & 4.0689E-04 & 0.9640 & 2.7367E-04 & 0.9969 & 2.7175E-04 & 0.9987 & 2.8275E-04 & 0.9991 \\
\hline
1.5 & 64 & 1.6696E-03 & -- & 1.6967E-03 & -- & 1.6847E-03 & -- & 1.6839E-03 & -- \\
    & 128 & 8.8651E-04 & 0.9133 & 8.7906E-04 & 0.9487 & 8.5788E-04 & 0.9736 & 8.5169E-04 & 0.9834 \\
    & 256 & 4.5771E-04 & 0.9537 & 4.4768E-04 & 0.9735 & 4.3266E-04 & 0.9875 & 4.2827E-04 & 0.9918 \\
    & 512 & 2.3261E-04 & 0.9765 & 2.2606E-04 & 0.9858 & 2.1752E-04 & 0.9921 & 2.1498E-04 & 0.9943 \\
\hline
1.9 & 64 & 1.5498E-03 & -- & 1.5600E-03 & -- & 1.5723E-03 & -- & 1.5782E-03 & -- \\
    & 128 & 8.0342E-04 & 0.9479 & 8.0002E-04 & 0.9634 & 7.9196E-04 & 0.9894 & 7.8715E-04 & 1.0036 \\
    & 256 & 4.0968E-04 & 0.9717 & 4.0908E-04 & 0.9677 & 4.0733E-04 & 0.9592 & 4.0636E-04 & 0.9539 \\
    & 512 & 2.0755E-04 & 0.9810 & 2.0649E-04 & 0.9863 & 2.0549E-04 & 0.9871 & 2.0549E-04 & 0.9837 \\
\hline
1.99 & 64 & 1.5094E-03 & -- & 1.5108E-03 & -- & 1.5130E-03 & -- & 1.5140E-03 & -- \\
     & 128 & 7.9818E-04 & 0.9192 & 7.9805E-04 & 0.9208 & 7.9781E-04 & 0.9233 & 7.9779E-04 & 0.9243 \\
     & 256 & 4.0060E-04 & 0.9946 & 4.0067E-04 & 0.9941 & 4.0079E-04 & 0.9932 & 4.0089E-04 & 0.9928 \\
     & 512 & 2.0409E-04 & 0.9730 & 2.0405E-04 & 0.9735 & 2.0399E-04 & 0.9743 & 2.0400E-04 & 0.9746 \\ 	
			\hline
		\end{tabular}
		\label{tab8}
	\end{center}
\end{table}	
\begin{table}[H]\footnotesize\tabcolsep=2.0pt
	\begin{center}
		\caption{The maximum norm errors between NL-IES \eqref{eq2.3} and L-IES \eqref{eq2.4} for Example 2
			with $\tau = h$.}
		\centering
		\begin{tabular}{cccccccccc}
			\hline
			& & \multicolumn{2}{c}{$\lambda = 0$} & \multicolumn{2}{c}{$\lambda = 1$}
			& \multicolumn{2}{c}{$\lambda = 5$} & \multicolumn{2}{c}{$\lambda = 10$} \\
			[-2pt] \cmidrule(lr){3-4} \cmidrule(lr){5-6} \cmidrule(lr){7-8} \cmidrule(lr){9-10} \\ [-11pt]
			$\alpha$ & $N$ & $Err(\tau,h)$ & $Order1$ & $Err(\tau,h)$ & $Order1$ & $Err(\tau,h)$ & $Order1$
			& $Err(\tau,h)$ & $Order1$ \\
			\hline
1.1 & 64 & 2.6283E-03 & -- & 1.8233E-03 & -- & 1.7754E-03 & -- & 2.0306E-03 & -- \\
    & 128 & 1.4530E-03 & 0.8551 & 9.7882E-04 & 0.8974 & 9.5216E-04 & 0.8989 & 9.4495E-04 & 1.1036 \\
    & 256 & 7.7570E-04 & 0.9055 & 5.1535E-04 & 0.9255 & 5.0437E-04 & 0.9167 & 5.0857E-04 & 0.8938 \\
    & 512 & 4.0344E-04 & 0.9431 & 2.6755E-04 & 0.9457 & 2.6380E-04 & 0.9350 & 2.7052E-04 & 0.9107 \\
\hline
1.5 & 64 & 1.6456E-03 & -- & 1.6295E-03 & -- & 1.6426E-03 & -- & 1.6259E-03 & -- \\
    & 128 & 8.7865E-04 & 0.9053 & 8.6650E-04 & 0.9112 & 8.4421E-04 & 0.9603 & 8.3915E-04 & 0.9542 \\
    & 256 & 4.5527E-04 & 0.9486 & 4.4503E-04 & 0.9613 & 4.3027E-04 & 0.9724 & 4.2567E-04 & 0.9792 \\
    & 512 & 2.3213E-04 & 0.9718 & 2.2560E-04 & 0.9801 & 2.1709E-04 & 0.9869 & 2.1452E-04 & 0.9886 \\
\hline
1.9 & 64 & 1.5331E-03 & -- & 1.5425E-03 & -- & 1.5498E-03 & -- & 1.5465E-03 & -- \\
    & 128 & 8.0010E-04 & 0.9382 & 7.9730E-04 & 0.9521 & 7.9150E-04 & 0.9694 & 7.8974E-04 & 0.9696 \\
    & 256 & 4.0891E-04 & 0.9684 & 4.0829E-04 & 0.9655 & 4.0667E-04 & 0.9607 & 4.0596E-04 & 0.9600 \\
    & 512 & 2.0742E-04 & 0.9792 & 2.0636E-04 & 0.9844 & 2.0531E-04 & 0.9861 & 2.0534E-04 & 0.9833 \\
\hline
1.99 & 64 & 1.4936E-03 & -- & 1.4938E-03 & -- & 1.4887E-03 & -- & 1.4758E-03 & -- \\
     & 128 & 7.9506E-04 & 0.9097 & 7.9510E-04 & 0.9098 & 7.9545E-04 & 0.9042 & 7.9606E-04 & 0.8905 \\
     & 256 & 3.9970E-04 & 0.9921 & 3.9974E-04 & 0.9921 & 3.9979E-04 & 0.9925 & 3.9976E-04 & 0.9937 \\
     & 512 & 2.0396E-04 & 0.9706 & 2.0392E-04 & 0.9711 & 2.0387E-04 & 0.9716 & 2.0388E-04 & 0.9714 \\
			\hline
		\end{tabular}
		\label{tab9}
	\end{center}
\end{table}

\begin{table}[H]\footnotesize\tabcolsep=3.0pt
	\begin{center}
		\caption{Results of different methods when $M = N$ for Example 2.}
		\centering
		\begin{tabular}{cccccccc}
			\hline
			& & \multicolumn{2}{c}{\rm{BS}} & \multicolumn{2}{c}{$\mathcal{I}$} & \multicolumn{2}{c}{$\mathcal{P}$} \\
			[-2pt] \cmidrule(lr){3-4} \cmidrule(lr){5-6} \cmidrule(lr){7-8} \\ [-11pt]
			$(\alpha, \lambda)$ & $N$ & $\mathrm{Iter1}$ & $\textrm{Time}$ & $(\mathrm{Iter1}, \mathrm{Iter2})$ & $\textrm{Time}$
			& $(\mathrm{Iter1}, \mathrm{Iter2})$ & $\textrm{Time}$ \\
			\hline
			(1.1, 0) & 129 & 5.0 & 0.539 & (5.0, 207.0) & 6.609 & (5.0, 3.8) & 0.167 \\
			& 257 & 5.0 & 4.022 & (5.0, 436.2) & 42.799 & (5.0, 6.0) & 0.815 \\
			& 513 & 5.0 & 37.223 & \ddag & \ddag & (5.0, 9.2) & 4.360 \\
			& 1025 & \dag & \dag & \ddag & \ddag & (5.0, 15.4) & 24.702 \\
			\hline
			(1.5, 0) & 129 & 4.0 & 0.471 & (4.0, 426.8) & 10.178 & (4.0, 11.3) & 0.331 \\
			& 257 & 4.0 & 3.413 & \ddag & \ddag & (4.0, 24.3) & 2.394 \\
			& 513 & 4.0 & 30.196 & \ddag & \ddag & (4.0, 62.0) & 22.520 \\
			& 1025 & \dag & \dag & \ddag & \ddag & (4.0, 154.8) & 190.537 \\
			\hline
			(1.9, 0) & 129 & 4.0 & 0.462 & (4.0, 880.5) & 22.810 & (4.0, 8.0) & 0.309 \\
			& 257 & 4.0 & 3.426 & \ddag & \ddag & (4.0, 21.8) & 3.258 \\
			& 513 & 4.0 & 30.240 & \ddag & \ddag & (4.0, 66.3) & 31.259 \\
			& 1025 & \dag & \dag & \ddag & \ddag & (4.0, 212.5) & 362.669 \\
			\hline
			(1.1, 5) & 129 & 5.0 & 0.553 & (5.0, 271.8) & 8.562 & (5.0, 2.8) & 0.122 \\
			& 257 & 5.0 & 4.105 & (5.0, 584.2) & 60.355 & (5.0, 4.2) & 0.536 \\
			& 513 & 5.0 & 36.435 & \ddag & \ddag & (5.0, 7.6) & 3.632 \\
			& 1025 & \dag & \dag & \ddag & \ddag & (5.0, 14.2) & 22.479 \\
			\hline
			(1.5, 5) & 129 & 4.0 & 0.470 & (4.0, 542.5) & 13.903 & (4.0, 6.5) & 0.189 \\
			& 257 & 4.0 & 3.460 & \ddag & \ddag & (4.0, 17.5) & 1.704 \\
			& 513 & 4.0 & 30.523 & \ddag & \ddag & (5.0, 46.4) & 20.737 \\
			& 1025 & \dag & \dag & \ddag & \ddag & (5.0, 136.2) & 206.973 \\
			\hline
			(1.9, 5) & 129 & 4.0 & 0.452 & (4.0, 907.0) & 24.141 & (4.0, 4.5) & 0.167 \\
			& 257 & 4.0 & 3.437 & \ddag & \ddag & (4.0, 15.5) & 2.241 \\
			& 513 & 4.0 & 30.595 & \ddag & \ddag & (4.0, 54.3) & 25.433 \\
			& 1025 & \dag & \dag & \ddag & \ddag & (4.0, 197.5) & 338.909 \\
			\hline
			(1.1, 10) & 129 & 5.0 & 0.568 & (5.0, 277.0) & 8.515 & (5.0, 2.8) & 0.109 \\
			& 257 & 5.0 & 4.069 & (5.0, 598.8) & 58.637 & (5.0, 3.6) & 0.409 \\
			& 513 & 5.0 & 36.612 & \ddag & \ddag & (5.0, 6.6) & 3.083 \\
			& 1025 & \dag & \dag & \ddag & \ddag & (6.0, 12.3) & 23.298 \\
			\hline
			(1.5, 10) & 129 & 5.0 & 0.549 & (5.0, 571.2) & 17.305 & (5.0, 4.2) & 0.158 \\
			& 257 & 5.0 & 4.097 & \ddag & \ddag & (5.0, 12.2) & 1.524 \\
			& 513 & 5.0 & 36.738 & \ddag & \ddag & (5.0, 38.6) & 17.222 \\
			& 1025 & \dag & \dag & \ddag & \ddag & (5.0, 124.0) & 191.098 \\
			\hline
			(1.9, 10) & 129 & 4.0 & 0.459 & (5.0, 841.2) & 26.614 & (4.0, 3.3) & 0.126 \\
			& 257 & 4.0 & 3.434 & \ddag & \ddag & (4.0, 11.5) & 1.689 \\
			& 513 & 4.0 & 30.590 & \ddag & \ddag & (4.0, 46.0) & 21.659 \\
			& 1025 & \dag & \dag & \ddag & \ddag & (4.0, 179.0) & 307.717 \\		
			\hline
		\end{tabular}
		\label{tab10}
	\end{center}
\end{table}

\section{Concluding remarks}
\label{sec6}\quad
The nonlinear all-at-once system arising from the nonlinear tempered fractional diffusion equations is studied.
Firstly, the two implicit schemes (i.e., NL-IES \eqref{eq2.3} and L-IES \eqref{eq2.4})
in Section \ref{sec2} are obtained through applying the finite difference method.
Then, the stabilities and first-order convergences of such schemes are analyzed in Section \ref{sec3}
under several suitable assumptions.
Secondly, for solving all the time steps in Eq. \eqref{eq2.3} simultaneously,
the nonlinear all-at-once system \eqref{eq2.5} is derived from it.
Then, Newton's method is employed to solve this system \eqref{eq2.5}.
Once the method is used to solve such the nonlinear system, the following two basic problems need to be addressed:
1. How to find a good initial value for Newton's method?
2. How to fast solving the Jacobian equations in Newton's method?
As for the first problem, the value, which is constructed by interpolating the solution
of L-IES (\ref{eq2.4}) on the coarse grid, is chosen as the initial guess.
For the second problem, the PBiCGSTAB method
with the preconditioner $P_{\ell}$ is employed to accelerate solving the Jacobian equations,
which is discussed in Section \ref{sec4}.
On the one hand, numerical examples in Section \ref{sec5} show that
the convergence orders of two proposed schemes (both in continuous and discontinuous cases)
can indeed reach $1$ in both time and space.
On the other hand, they also indicate that our preconditioning strategy
is effective for solving \eqref{eq2.5} with continuous or discontinuous coefficients.
However, the performance of $P_{\ell}$ are not satisfactory.
The reason may be that the diagonal block matrix in the Jacobian matrix
(can be rewritten as $\mathcal{A}$ plus this diagonal block matrix),
which is resulted from $-\tau \frac{\partial \bm{f(u)}}{\partial \bm{u}}$,
is not considered when designing $P_{\ell}$ in this work.
Thus, a preconditioner designed with considering such the diagonal block matrix
may be more effective to solve Eq. \eqref{eq2.5}.
In the future work, we will study along with this direction and give some relative theoretical analysis.
\section*{Acknowledgments}
\quad
\textit{This research is supported by the National Natural Science Foundation of China (Nos. 61772003, 61876203 and 11801463)
and the Fundamental Research Funds for the Central Universities (No. ZYGX2016J132).}
%
\bibliographystyle{elsarticle-num}

\end{document}